\documentclass[11pt]{amsart}

\newtheorem{thm}{Theorem}[section]

\newtheorem{assumption}[thm]{Assumption}

\usepackage{amsmath}
\usepackage{amsfonts}
\usepackage{amssymb}
\usepackage{amsthm}
\usepackage{mathtools} 
\usepackage{latexsym}

\usepackage{enumerate}
\usepackage{enumitem}
\usepackage{cancel}
\usepackage{cases}
\usepackage{empheq}
\usepackage{multicol}

\usepackage{float}

\usepackage{color}

\usepackage{mathrsfs}
\usepackage{fontenc} 
\usepackage{inputenc} 
\usepackage{verbatim}

\theoremstyle{plain}
\newtheorem{theorem}{Theorem}[section]

\newtheorem{lemma}[theorem]{Lemma}

\theoremstyle{definition}
\newtheorem{definition}[theorem]{Definition}

\theoremstyle{remark}
\newtheorem{remark}[theorem]{Remark}

\numberwithin{equation}{section} %
\numberwithin{figure}{section}   %

\usepackage[square,comma,numbers,sort]{natbib}
\usepackage[colorlinks=true, pdfborder={ 0 0 0}]{hyperref}
\hypersetup{urlcolor=blue, citecolor=red}
\usepackage{url}

\newcommand{\field}[1]{\mathbb{#1}}

\newcommand{\nR}{\field{R}}

\newcommand{\tu}{\tilde{u}}

\newcommand{\tb}{\tilde{b}}
\newcommand{\ttheta}{\tilde{\theta}}

\newcommand{\tj}{ \abs{\grad\times B} - \abs{\grad\times b}  }

\newcommand{\ej}{\abs{\grad\times b}}

\newcommand{\nT}{\mathbb T}

\newcommand{\lp}{\left(}
\newcommand{\rp}{\right)}

\newcommand{\la}{\left\langle}
\newcommand{\ra}{\right\rangle}

\newcommand{\dt}{\partial_t}
\newcommand{\lap}{\triangle}
\newcommand{\grad}{\nabla}
\newcommand{\bze}{\zeta^\epsilon}

\newcommand{\T}{\mathbb{T}^2}
\newcommand{\jj}{\abs{\grad\times B}}

\newcommand{\abs}[1]{\left\lvert#1\right\rvert}

\newcommand{\ip}[2]{\left<#1,#2\right>}

\newcommand{\normLp}[2]{\|#2 \|_{L^{#1}}}
\newcommand{\normHs}[2]{\|#2 \|_{H^{#1}}}

\newcommand{\normL}[1]{\|#1\|_{L^2}}
\newcommand{\normH}[1]{\|#1\|_{H^1}}
\newcommand{\fournorm}[1]{\|#1\|_{L^4}}
\newcommand{\infnorm}[1]{\|#1\|_{L^{\infty}}}
\newcommand{\Htwonorm}[1]{\|#1\|_{H^{2}}}

\renewcommand{\mu}{\eta}

\usepackage{url}

\usepackage{placeins} 

\newcounter{my_counter}
\title[Calmed Ohmic Heating]{
Calmed Ohmic Heating for the 2D Magnetohydrodynamic-Boussinesq System: Global Well-posedness and Convergence}

\date{\today}

\author{Matthew Enlow}
\address[Matthew Enlow]{Department of Mathematics, 
                University of Nebraska--Lincoln,
        Lincoln, NE 68588-0130, USA}
\email[Matthew Enlow]{menlow2@huskers.unl.edu}

\author{Adam Larios}
\address[Adam Larios]{Department of Mathematics, 
                University of Nebraska--Lincoln,
        Lincoln, NE 68588-0130, USA}
\email[Adam Larios]{alarios@unl.edu}

\author{Yuan Pei}
\address[Yuan Pei]{Department of Mathematics, 
                Western Washington University
        Bellingham, WA 98225, USA}
\email[Yuan Pei]{peiy@wwu.edu}


\begin{document}
\begin{abstract}
When an electric current runs through a fluid, it generates heat via a process known as ``Ohmic heating'' or ``Joule heating.''  While this phenomenon, and its quantification known as Joule's Law, is the first studied example of heat generation via an electric field, many difficulties still remain in understanding its consequences.  In particular, a magnetic fluid naturally generates an electric field via Amp\`ere's law, which heats the fluid via Joule's law.  This heat in turn gives rise to convective effects in the fluid, creating complicated dynamical behavior.  This has been modeled (in other works) by including an Ohmic heating term in the Magnetohydrodynamic-Boussinessq (MHD-B) equation.  However, the structure of this term causes major analytical difficulties, and basic questions of well-posedness remain open problems, even in the two-dimensional case.  Moreover, standard approaches to finding a globally well-posed approximate model, such as filtering or adding high-order diffusion, are not enough to handle the Ohmic heating term. In this work, we present a different approach that we call ``calming'', which reduces the effective algebraic degree of the Ohmic heating term in a controlled manner.  We show that this new model is globally well-posed, and moreover, its solutions converge to solutions of the MHD-B system with the Ohmic heating term (assuming that solutions to the original equation exist), making it the first globally well-posed approximate model for the MHD-B equation with Ohmic heating.
\end{abstract}

\maketitle
\thispagestyle{empty}

\section{Introduction}\label{secInt}
\noindent
Turbulent fluids that involve multi-physical processes, such as magnetic effects and/or heat conduction are a fascinating and challenging area of research. One way to understand such fluids is by studying the magnetohydrodynamic (MHD) partial differential equations with an additional thermal equation coupled to the original system \cite{Larios_Pei_2017_MHDB_PS}.  Such equations have been considered in many works, often under the name ``MHD-Boussinesq'' (MHD-B) equations \cite{Pan_2022_MHDB_axisym,Liu_Deng_Lin_Sun_2022_3D_Bouss_MHD, wang2023global, ye2021global, Liu_Bian_Pu_2019global, Bian_Gui_2016_MHDBouss, Chamorro_He_2021_partial_MHD, liu2023asymptotic},
but the heating effect of the electrical current  (so-called ``Ohmic heating'' or ``Joule heating'') is neglected, save for in a small number of papers \cite{Manley_Treve_1981, chen2004combined, goud2021ohmic, sibanda2010steady, soid2017unsteady, pal2010buoyancy, su2012mhd}.
Therefore, we are interested in the case of the MHD-B system \textit{with} Ohmic heating (MHD-B$\Omega$).  We note that \textit{even in the two-dimensional (2D) case, global well-posedness for this system is a completely open problem}.  
Thus, in the present work, we consider the 2D MHD-B system with a \textit{modified} Ohmic heating term.  We prove that our modified system is globally well-posed, and also that, at least before any potential blow-up time of the original system, solutions of our modified system converge to solutions of the MHD-B$\Omega$ system, assuming solutions to the latter equation exist.  

The difficulty with the MHD-B$\Omega$ system is that the Ohmic heating term is already quadratic and non-negative, destroying any hope of closing energy estimates, or of finding any conserved quantity.  Moreover, there is no natural scaling to the equation.  Indeed, the Ohmic term $|\nabla\times b|^2$ scales exactly like the well-known Hall term $\nabla\times((\nabla\times b)\times b)$ of the Hall-MHD system, 
for which global well-posedness is also a highly challenging open problem in 2D
\footnote{Technically, the Hall-MHD system only makes sense in 3D, or in the so-called ``two-and-a -half dimensional'' (2.5D) case, where the spatial inputs are 2D, but the output dimension is 3D Global well-posedness for the 2.5D Hall-MHD is the open problem which we are referring to here.}.  See, for example, \cite{dai2016_HallMHD_reg, dai2019_hallMHD_longtime, dai2020_hallMHD_lwp, dai2021_hallMHD_nonuniqueLH, dai2022_HallMHD_gwp} for a background of the wellposedness and regularity results obtained for the Hall-MHD system in dimensions $n \geq 2$ in the deterministic case, and \cite{yamazaki2017_hallMHD_stochastic, yamazaki2019_hallMHD_markov, yamazaki2019_hallMHD_wpLevy, yamazaki2019_hallMHD_stochasticRemarks, yamazaki2020_hallMHD_irreducibility, yamazaki2022_hallMHD_ergodicity, yamazaki_rahman2022_hallMHDremarks} for results on the stochastic Hall-MHD system.
Indeed, the MHD-B$\Omega$ system appears to be even more challenging than the Hall-MHD system: the Hall term vanishes in $L^2$-energy estimates allowing for a conserved quantity, but no analogous situation is known to hold for the Ohmic term.

The above difficulties with the Ohmic term mean that standard modifications, such as adding hyper diffusion, filtering the advective term, etc., have little chance of allowing for a proof of global well-posedness, since these techniques merely control the growth of gradients (i.e., the growth of small scales).  What is needed is some way to control the Ricatti-like\footnote{Recall that the Ricatti-type equation $\frac{dy}{dt}=y^{1+\epsilon}$ blows up in finite time for any $\epsilon>0$ and positive initial data.} nature of the system.  Note that similar issues arise in controlling the vorticity in the 3D Navier-Stokes equations (NSE), and controlling the solution in the 2D Kuramoto-Sivashinsky equation (KSE).  Hence, we employ a technique that we call ``algebraic calming,''  developed in 
\cite{Enlow_Larios_Wu_2023_KSE} in the context of the 2D KSE.  However, rather than apply the calming function to the advective term, as is done in \cite{Enlow_Larios_Wu_2023_KSE}, the innovation of the present work is to apply the calming function to the Ohmic term.  We describe this in detail below.

Consider the domain $\mathbb{T}^2$ the two-dimensional periodic space $\mathbb{R}^2/\mathbb{Z}^2 = [0, 1]^2$. 
For $T>0$, the 2D MHD-Boussinesq system with full fluid viscosity $\nu> 0$, magnetic resistivity $\mu> 0$, and thermal diffusion $\kappa > 0$
over $\nT^2\times [0, T)$, is given by
\begin{equation}
       \left\{
       \begin{aligned}
       &\dt u 
       - \nu\Delta u 
       + (u\cdot\nabla) u 
       + \nabla p 
       = 
       (b\cdot\nabla) b
       +
       g\theta \vec{e}_2,
     \\&
       \dt b 
       - \mu\Delta b 
       + (u\cdot\nabla) b
       = 
       (b\cdot\nabla) u,
     \\&
       \dt \theta 
       -
        \kappa\Delta\theta
       + 
       (u\cdot\nabla) \theta  
       = 
       \alpha\mu \bze(|\nabla\times b|)|\nabla\times b|, 
    \\&
    \nabla \cdot u = 0 = \nabla \cdot b,
    \end{aligned}
    \right.
    \label{Sys1}
  \end{equation}
where the constant $g>0$ has unit of force, and is proportional to the constant of gravitational acceleration.  
We denote $x = (x_1, x_2)$, and $\vec{e}_2$ to be the unit vector in the $x_2$ direction, i.e., $\vec{e}_2 = (0, 1)^{T}$.
Here and henceforth, $u = u(x, t) = (u_1(x, t), u_2(x,t))$ is the unknown velocity field of a viscous incompressible fluid, 
with divergence-free initial data $u(x, 0)=u_0$;
$b = b(x, t) = (b_1(x, t), b_2(x,t))$ is the unknown magnetic field, 
with divergence-free initial data $b(x, 0)=b_0 $;
and the scalar $p=p(x, t)$ represents the unknown pressure, while $\theta = \theta(x, t)$ can be thought of as the unknown temperature fluctuation, with initial value $\theta_{0} = \theta(x, 0)$.

Here, we assume that the calming function $\bze$ is a bounded, smooth, globally Lipschitz function with Lipschitz constant $L_\epsilon >0$ which also satisfies 
\begin{align}\label{calm_bounded}
    |\bze(x)| + |\nabla \bze(x)|\leq M_{\epsilon}
\end{align}
for some $M_\epsilon > 0$ such that $\lim_{\epsilon \to 0^+} M_\epsilon = \infty$.
Also, we assume the calming function $\bze$ approximates the identity function in the following sense:
\begin{align}
    \label{calm_convergence_rate}
    |x - \bze(x)| \leq C \epsilon^\gamma |x|^\beta
\end{align}
For some $C > 0$, $\gamma > 0$ and for $\beta \geq 1$. We take $\bze$ to be any function which satisfies these properties, though we give several concrete examples below.
\begin{minipage}{0.48\textwidth}
\begin{align} \label{eta_choices}
\bze(x) = 
\begin{cases}
  \bze_1(x) &:= \cfrac{x}{1+\epsilon|x|},\\
  \bze_2(x) &:= \cfrac{x}{1+\epsilon^2|x|^2},\\
  \bze_3(x) &:= \cfrac{1}{\epsilon}\arctan(\epsilon x),\\
  \bze_4(x) &:= q^\epsilon(\abs{x}) \frac{x}{\abs{x}},
\end{cases}
\end{align}
\end{minipage}

\vspace{0.2 in}
where $q^\epsilon$ is a piecewise function given by 
\begin{align}\notag
q^\epsilon(r) = 
    \begin{cases}
        r, \quad  
        & 0 \leq  r < \frac{1}{\epsilon}, \\ 
        -\frac{\epsilon}{2}\lp r - \frac{2}{\epsilon} \rp^2 + 
        \frac{3}{2\epsilon}, \quad  
        & \frac{1}{\epsilon} \leq r < \frac{2}{\epsilon}, \\ 
        \frac{3}{2\epsilon} , \quad 
        & r \geq  \frac{2}{\epsilon}. 
    \end{cases}
\end{align}

\begin{remark}
    One might expect that the modification $\bze(x) = \abs{x}^{1-\epsilon}$ would suffice to show wellposedness of this system. However, this modified system is only wellposed for $\epsilon \geq 1$. Our calmed Ohmic heating system has the benefit of being wellposed for all $\epsilon > 0$ and converges to the original system as $\epsilon \to 0^+$.
\end{remark}

With the same parameters, we now give below the original two-dimensional Boussinesq-MHD equations with Ohmic heating effect but without the calming modification.
\begin{equation}
       \left\{
       \begin{aligned}
       &\dt U 
       - \nu\Delta U 
       + (U\cdot\nabla) U 
       + \nabla P 
       = 
       (B\cdot \nabla) B
       +
       g\Theta \vec{e}_2,
     \\&
       \dt B 
       - \mu\Delta B 
       + (U\cdot\nabla) B
       = 
       (B\cdot\nabla) U,
     \\&
       \dt \Theta 
       -
        \kappa\Delta\Theta
       + 
       (U\cdot \nabla) \Theta  
       = 
       \alpha\mu |\nabla\times B|^2, 
    \\&
    \nabla \cdot U = 0 = \nabla \cdot B,
    \end{aligned}
    \right.
    \label{Sys2}
  \end{equation}

  \subsection{A Toy Model}
To shed light on some of the ideas in the present work, we consider the following ``calmed'' version of the well-known Ricatti equation $y'=y^2$; namely $y'=y\bze(y)$, where for the sake of algebraic simplicity, we consider the calming function $\bze=\bze_2$ in \eqref{eta_choices}. That is, for $\epsilon>0$ and $y_0>0$, we consider the initial value problem
\begin{align}
    \begin{cases}
        \,\,\, \frac{dy}{dt} \!\!\!\!&= y\,\frac{y}{1+\epsilon^2y^2},\\
        y(0)\!\!\!\!&=y_0.
    \end{cases}
\end{align}
The unique solution to this calmed equation (using $y_0>0$) is given by
  \begin{align}
      y_\epsilon(t) 
      = 
      \frac{\epsilon^2y_0^2+ty_0 - 1+ \sqrt{(\epsilon^2y_0^2+ty_0 - 1)^2+4\epsilon^2y_0^2}}{2 \epsilon^2y_0}, 
  \end{align}
which exists globally; i.e., for all $t\in\nR$.  For fixed $t\in[0,\frac{1}{y_0})$, $\lim_{\epsilon\rightarrow0}y_\epsilon(t) = \frac{y_0}{1-ty_0}$, so the usual solution of the Ricatti equation is attained.  (In fact, the convergence is uniform on $(-\infty,T]$ for every $T<\frac{1}{y_0}$).  
However, for $t\geq\frac{1}{y_0}$, the limit diverges to $+\infty$.  To summarize, the calmed solution is global, and it gives the correct limiting solution of the original equation on the interval of existence and uniqueness of the original equation, but the limit does not provide a notion of solution to the original equation beyond this interval, at least in this example.

\section{Preliminaries and Main Results}\label{pre_main}

In this section, we first provide the essential preliminaries that we need; then, we list the main theorems and the organization of this paper.

\begin{remark}
Note that all of the results in the present work also hold easily \textit{mutatis mutandis} for the so-called ``two-and-a -half dimensional'' case, that is the case where $x=(x_1,x_2)$ is still two-dimensional, but the the outputs are three-dimensional, i.e., 
$$u = u(x, t) = (u_1(x, t), u_2(x,t),u_3(x,t))$$ 
and 
$$b = b(x, t) = (b_1(x, t), b_2(x,t),b_3(x,t)).$$  
This is because the only role dimenionality plays in our analysis is in Sobolev estimates, which depend only on the input dimension.  For the sake of simplicty, we present only the 2D case.
\end{remark}

\subsection{Preliminaries}\label{pre}

Throughout this paper, we denote by $f_{x_1}$ the partial derivative of $f(x_1,x_2)$ with respect to $x_1$, and likewise for $g_{x_2}$, etc.. We also denote the Lebesgue and the Sobolev spaces by $L^p$ for $0\leq p \leq \infty$ and $H^s = W^{s,2}$ with $s>0$, respectively. Let $\mathcal V$ be the set of all $L$-periodic trigonometric polynomials from $\nR^2$ to $\nR^2$ that are divergence free with zero average. We denote by $H$, $V$, and $V'$ the closures of $\mathcal V$ in the $L^2(\nT^2)^2$ and $H^1(\nT^2)^2$ norms, and the dual space of $V$, respectively, with inner products on $H$ and $V$ as
$$(u, v) = \sum_{i=1}^2\int_{\mathbb{T}^2}u_{i}v_{i}\,dx \text{ \,\,and\,\,  } (\nabla u, \nabla v) = \sum_{i, j=1}^2\int_{\mathbb{T}^2}\partial_{j}u_{i}\partial_{j}v_{i}\,dx,$$
respectively, associated with the norms 
$$\| u \|_{H}=(u, u)^{1/2} \,\,\text{and}\,\, \| u \|_{V}=(\nabla u, \nabla u)^{1/2},$$ where 
$$u(x) = (u_1(x_1, x_2), u_2(x_1, x_2))\,\,\text{and}\,\,v(x) = (v_1(x_1, x_2), v_2(x_1, x_2)).$$ 
For the sake of convenience, we use $\normL{u}$ and $\normH{u}$ to denote the above norms in $H$ and $V$, respectively. 

We define the Stokes operator $A:= -P_{\sigma}\Delta$ 
with domain $\mathcal{D}(A):= H^2\cap V$, 
where $P_{\sigma}$ is the orthogonal Leray-Helmholtz projector from $L^2$ to $H$. 
Notice that on our domain $\Omega = \mathbb{T}^2$, 
we have $A = -\Delta P_{\sigma}$.
Moreover, the Stokes operator can be extended 
to a linear operator from $V$ to $V'$ as 
$$\left<Au, v\right> = (\nabla u, \nabla v) \text{  for all  } v\in V.$$
Moreover, we have the well-known fact that $A^{-1} : H \hookrightarrow \mathcal{D}(A)$ 
is a positive-definite, self-adjoint, and compact operator from $H$ into itself,
thus, $A^{-1}$ possesses an orthonormal basis of positive eigenfunctions $\{ w_{k}\}_{k=1}^{\infty}$ in $H$, corresponding to a sequence of non-increasing sequence of eigenvalues. 
Therefore, $A$ has non-decreasing eigenvalues $\lambda_{k}$, 
i.e., $0 < \lambda_1 \leq \lambda_2, \dots$. 
Then, we have the following Poincar\'e inqualities:
\begin{align}\label{Poincare}
    \lambda_1\normL{u}^2
    \leq
    \normL{\grad u}^2 \text{\,\, for\,\,} u\in V;\quad
    \lambda_1\normL{\grad u}^2
    \leq
    \normL{Au}^2 \text{\,\, for\,\,} u\in \mathcal{D}(A).
\end{align}
Thus, $\normL{\grad u}$ is equivalent to $\normH{u}$ for mean-free functions (in particular, when $u$ is a derivative, thanks to the periodic boundary conditions). 

We frequently use the following Agmon's inequality in two-dimensional space.
\begin{equation}\label{agmon}
    \infnorm{u}^2 \leq c\normL{u} \Htwonorm{u},
\end{equation}
for all $u\in H^2$. 

Also, we state the following Ladyzhenskaya inequality 
\begin{equation}\label{ladyzhenskaya}
    \fournorm{u}^2 \leq c\normL{u} \normH{u},
\end{equation}
for all $u\in V$, which is a special case of the following interpolation result 
that is frequently used in this paper (see, e.g., \cite{Nirenberg_1959_AnnPisa} for a detailed proof).
Assume $1 \leq q, r \leq \infty$, and $0<\gamma<1$.  
For $v\in L^q(\mathbb{T}^{n})$, such that  $\partial^\alpha v\in L^{r} (\mathbb{T}^{n})$, for $|\alpha|=m$, then 
\begin{align}\label{PT1}
\|\partial_{s}v\|_{L^{p}} \leq C\|\partial^{\alpha}v\|_{L^{r}}^{\gamma}\| v\|_{L^{q}}^{1-\gamma},
\,\,\,\text{where}\,\,\,
\frac{1}{p} - \frac{s}{n} = \left(\frac{1}{r} - \frac{m}{n}\right) \gamma+ \frac{1}{q}(1-\gamma).
\end{align}

A special case of the above Sobolev-type inequalities, frequently used in this paper, is summarized the following lemma. 
\begin{lemma}
  For $f, g, h\in H^{1}(\nT^2)$, we have 
  $$\iint_{\nT^2}|fgh|\,dx_1\,dx_2 \leq C \normL{f}\normL{g}^{1/2}\normH{g}^{1/2}\normL{h}^{1/2}\normH{h}^{1/2},$$
  for some constant $C>0$. 
\end{lemma}
Additionally, we use the following result which can be found in, e.g., \cite[Problem 5.10.17]{Evans_2010}: 
\begin{lemma}\label{evans_lem}
    If $u \in W^{1,p}(U)$ on a bounded domain $U$ for $1<p<\infty$, then $\abs{u} \in W^{1,p}(U)$.
\end{lemma}

We use the following standard notation for the projection of the nonlinear term,
\begin{align*}
     \mathcal{B}(u, v) := P_{\sigma}((u\cdot\nabla)v),\quad (u, v \in \mathcal{V}),
\end{align*}
which can be extended to a continuous bilinear map
$\mathcal{B} : V \times V \to V'$. Moreover, the following lemma holds (see, e.g., \cite{Constantin_Foias_1988,Temam_2001_Th_Num,Foias_Manley_Rosa_Temam_2001} for details).
\begin{lemma}
For all $u,v,w\in V$, 
    \begin{subequations}
      \begin{align}
      \label{symm1}
      \ip{\mathcal{B}(u,v)}{w}_{V'} &= -\ip{\mathcal{B}(u,w)}{v}_{V'}, \\
      \label{symm2}
      \ip{\mathcal{B}(u,v)}{v}_{V'} &= 0.
    \end{align}
  \end{subequations}
Also, for all $u,v,$ and $w$ in the largest spaces $H$, $V$, or $\mathcal{D}(A)$, for which the right-hand sides of the inequalities below are finite,
  \begin{subequations}
    \begin{align}
      \label{B:424}
      |\ip{\mathcal{B}(u,v)}{w}_{V'}|
      &\leq 
      C\normL{u}^{1/2} \normH{u}^{1/2} \normH{v} \normL{w}^{1/2} \normH{w}^{1/2},
      \\
      \label{B:442}
      |\ip{\mathcal{B}(u,v)}{w}_{V'}|
      &\leq 
      C\normL{u}^{1/2} \normH{u}^{1/2} \normH{v}^{1/2} \normL{Av}^{1/2} \normL{w},    
      \\
      \label{B:inf22}
      |\ip{\mathcal{B}(u,v)}{w}_{V'}|
      &\leq 
      C\normL{u}^{1/2} \normL{Au}^{1/2} \normH{v} \normL{w}.
    \end{align}
Moreover, due to the periodic boundary conditions, it holds (one in two-dimensions) that
  \end{subequations}
  \begin{align}\label{enstrophy_miracle}
    \ip{\mathcal{B}(w, w)}{Aw} = 0,\; w\in \mathcal{D}(A),
\end{align}
  and the following Jacobi identity holds 
  \begin{align}\label{jacobi}
    \ip{\mathcal{B}(u, w)}{Aw} + \ip{\mathcal{B}(w, u)}{Aw} +\ip{\mathcal{B}(w, w)}{Au} = 0.
  \end{align}
\end{lemma}

For the sake of completeness, we state the following uniform Gr\"{o}nwall's inequality, proved in \cite{Jones_Titi_1992} (see also \cite{Farhat_Lunasin_Titi_2017_Horizontal} and the references therein), which will be used frequently throughout the paper. 
\begin{lemma}
\label{Gronwall}
Suppose that $Y(t)$ is a locally integrable and absolutely continuous function that satisfies the following: 
$$\frac{d Y}{d t} + \alpha(t) Y \leq \beta(t), \quad\text{ a.e. on } (0, \infty), $$
such that 
$$\liminf_{t \to \infty} \int_{t}^{t+\tau} \alpha(s)\,ds \geq \gamma, \quad\quad\quad \limsup_{t \to \infty} \int_{t}^{t+\tau} \alpha^{-}(s)\,ds < \infty, $$
and 
$$\lim_{t \to \infty} \int_{t}^{t+\tau} \beta^{+}(s)\,ds = 0, $$
for fixed $\tau > 0$, and $\gamma > 0$, 
where 
$\alpha^{-} = \max\{-\alpha, 0\}$ 
and 
$\beta^{+} = \max\{\beta, 0\}$. 
Then, $Y(t) \to 0$ at an exponential rate as $t \to \infty$.
\end{lemma}

To prove convergence, we will need to use the following bootstrap principle (see, e.g., \cite[p. 20]{Tao_2006}).
\begin{lemma} \label{boot}
	Let $T>0$. Assume that two statements $C(t)$ and $H(t)$ with $t\in [0,T]$ 
	satisfy the
	following conditions:
	\begin{enumerate}[label=(\alph*)]
		\item\label{lem_a} If $H(t)$ holds for some $t\in [0,T]$, then $C(t)$ holds for 
		the same $t$;
		\item\label{lem_b} If $C(t_0)$ holds for some $t_0\in [0,T]$, then $H(t_0)$ holds 
		for $t$ in a neighborhood of $t_0$;
		\item\label{lem_c} If $C(t_m)$ holds for $t_m\in [0,T]$ and $t_m \to t_c$, then $C(t_c)$ holds;
		\item\label{lem_d} $H(t)$ holds for at least one $t \in [0,T]$.
	\end{enumerate}
	Then $C(t)$ holds for all $t\in[0,T]$.
\end{lemma}

\subsection{Main Results}\label{main}

We state in this section our major theorems and the outline of the remaining of this paper. Throughout, we assume that the calming function $\bze$ is a Lipschitz function satisfying Conditions \eqref{calm_bounded} and \eqref{calm_convergence_rate} with fixed $\epsilon>0$.

First, we define our notions of a solution to system \eqref{Sys1}.

\begin{definition}\label{def1}
   For $T>0$ and initial data $(u_0, b_0, \theta_0) \in V\times V\times L^2$, we say that $(u,b,\theta)$ is a \textit{semi-strong solution} to the calmed system \eqref{Sys1} on $[0,T]$ if $(u,b,\theta)$ satisfies \eqref{Sys1} in a functional sense in $L^2(0,T; H)\times L^2(0,T; H) \times L^2(0,T; H^{-1})$,  if the initial data is satisfied in $C([0,T];V) \times C([0,T];V) \times C([0,T];L^2)$, and moreover,  
    $$u, b\in C([0,T]; V) \cap L^2((0,T);H^{2}\cap V),$$
    with 
    $$\dt u, \dt b \in L^2(0,T; H),$$
    and 
    $$\theta\in C([0,T]; L^2) \cap L^2((0,T); H^1),$$
    with 
    $$\dt \theta \in L^2(0,T; H^{-1}).$$
    
\end{definition}

\begin{remark}
    We use the nomenclature \textit{semi-strong} due to that fact that we only take $\theta_0 \in L^2$, typically associated with weak solutions, while $u,b$ are in the more regular space $V$, typically associated with strong solutions of the NSE or MHD equations. Moreover, it is worth noting that we are not able to provide a proof of uniqueness for semi-strong solutions, but we can prove uniqueness of strong solutions, as defined below.
\end{remark}

\begin{definition}\label{def2}
   For $T>0$ and initial data $(u_0, b_0, \theta_0) \in (H^2\cap V)\times (H^2\cap V)\times H^2$, we say that $(u,b,\theta)$ is a \textit{strong solution} to the calmed system \eqref{Sys1} on $[0,T]$ if $(u,b,\theta)$ 
   is a semi-strong solution to \eqref{Sys1} with initial data $(u_0, b_0, \theta_0)$ and in addition satisfies 
   $$u, b\in C([0,T]; H^2\cap V) \cap L^2(0,T; H^3\cap V), \quad 
   \theta \in L^\infty(0,T; H^2) \cap L^2(0,T;H^{3}),$$ where the initial data is satisfied in the sense of $C([0,T]; H^2\cap V)$ for $u_0$ and $b_0$, and in the sense of $C([0,T]; L^2)$ for $\theta_0$.

\end{definition}
\begin{remark}\label{rmk_reg_problem}
The regularity in Definition \ref{def2} implies $\nabla\times b\in C([0,T],V)\cap L^2(H^2\cap V)$, but this only implies $|\nabla\times b|\in C([0,T],H^1)$, which is not enough to infer that $\partial_t\theta\in L^2(0,T;H^1)$, and hence we are not able to apply the Aubin-Lions Lemma to obtain $\theta\in C(0,T;H^2)$.  Hence, we only ask that for strong solution the initial data $\theta_0$ be satisfied in the sense of $C([0,T]; L^2)$, a redundant condition from that of semi-strong solutions.
\end{remark}
Our first major result is the global existence of a semi-strong solution, although we are not able to prove uniqueness of these solutions (perhaps analogously with the fact that it is not know that weak Leray-Hopf solutions of the 3D Navier-Stokes are unique.  Similar remarks hold for the calmed Navier-Stokes equations \cite{Enlow_Larios_Wu_2023_NSE}).
\begin{theorem}\label{thm_globExi}
    Given arbitrary time $T>0$ and initial data $(u_0, b_0, \theta_0) \in V\times V\times L^2$, semi-strong solutions to the calmed system \eqref{Sys1} exist on $[0,T]$.
\end{theorem}

\begin{remark}
Estimates higher than order $s=2$ are not as straight-forward to obtain as they are in the NSE (or MHD-B equations), due to the presence of the calming term.  The reason for this is that, after integration by parts during estimates, higher-order derivatives of the calming term are complicated by the effect of the chain rule.  (To be clear, \textit{without} the calming term, \textit{no} estimates are known to be possible for MHD-B$\Omega$.)  It may be possible to handle this with, e.g., Schauder estimates, but such a discussion would likely distract from the main goal of the paper, which is to establish global well-posedness of an approximate model for Ohmic heating, and its convergence.  Hence, we focus on the case $s=2$.
\end{remark}

We also show $H^2$ regularity of solutions with more regular initial data.

\begin{theorem}\label{thm_hor}
    Given arbitrary time $T>0$, and initial conditions $u_0, b_0\in H^2\cap V$, $\theta_0\in H^2$, there exists a strong solution $(u,b,\theta)$ to System~\ref{Sys1} 
    which satisfies the conditions of Definitions \ref{def1} and \ref{def2}.

\end{theorem}
\begin{theorem}\label{thm_uni}
    Given arbitrary time $T>0$ and initial data $u_0, b_0, \theta_0 \in (H^2\cap V)\times (H^2 \cap V)\times L^2$, strong solutions to the calmed system \eqref{Sys1} are unique.
\end{theorem}

In order to obtain error estimates of the solution to System~\ref{Sys1} compared to that of System~\ref{Sys2}, we need the following local-in-time well-posedness assumption of the original system.  

\begin{assumption}[Short-time existence and uniqueness of strong solutions to the original MHD-Boussinesq System with Ohmic heating.]\label{thm_short_original}
    For $U_0$, $B_0\in H^2\cap V$ and $\Theta_0\in H^2$, there is a time $T^*>0$ 
    
    for which a unique solution $(U,B,\Theta)$ to System~\ref{Sys2} exists on $[0,T]$ for any $T\in(0,T^*)$, with 
    \begin{align*}
        U, B&\in C([0,T]; H^2\cap V) \cap L^2(0,T;H^{3}\cap V),
        \\
        \partial_tU, \partial_tB&\in  L^2(0,T;H),
    \end{align*}
    and 
    \begin{align*}
    \Theta&\in C([0,T]; H^2) \cap L^2(0,T;H^{3}),
    \\
    \partial_t\Theta&\in  L^2(0,T;L^2), 
    \end{align*}
    where the initial data is satisfied in the sense 
    of  $C([0,T];V) \times C([0,T];V) \times C([0,T];L^2)$.\end{assumption}
\begin{remark}
    Recall that, as pointed out in the introduction, the short-time existence of solutions to system \eqref{Sys2} is currently an open problem. Hence, our convergence results rest on the assumption that Assumption \ref{thm_short_original} holds, although this is not known.
\end{remark}

The next theorem concerns the convergence of the solution to System~\ref{Sys1} with calmed Ohmic heating to that of the original Boussinesq-MHD system without the calming mechanism (\ref{Sys2}), on the time-interval of existence of solutions of the latter. 
\begin{theorem}[Error analysis and convergence of the solution of (\ref{Sys1}) to that of (\ref{Sys2})]\label{thm_conv}
    Let $(U, B, \Theta)$ be the solution to system \eqref{Sys2} satisfying the conditions of Assumption \ref{thm_short_original} for $T > 0$ with initial data $U_0, B_0 \in H^2 \cap V$, $\Theta_0 \in H^2$.
    Assume that $\bze$ is a calming function which is Lipschitz and satisfies \eqref{calm_bounded} and \eqref{calm_convergence_rate}, and let $(u, b, \theta)$ be the strong solution to \eqref{Sys1}
    on $[0,T]$

    with initial data 
    $$u_0=U_0,\, b_0=B_0, \,\theta_0=\Theta_0.$$
    Then for all $t\in [0, T]$, 
    \begin{align}
        \normH{U(t)-u(t)}^2 +\normH{B(t)-b(t)}^2 +\normL{\Theta(t)-\theta(t)}^2 &\leq c_1\epsilon^{2\gamma}, \label{thm_conv_rate1} \\
        \int_0^T \normHs{2}{U(t)-u(t)}^2 +\normHs{2}{B(t)-b(t)}^2 +\normHs{1}{\Theta(t)-\theta(t)}^2 dt&\leq c_2\epsilon^{2\gamma}, \label{thm_conv_rate2} 
    \end{align}
    where the constants $c_1$ and $c_2$ are shown in \eqref{c_1} and \eqref{c_2} respectively and depend on $g,\nu,\mu,\kappa, \alpha$, $\normHs{2}{U}, \normHs{2}{B}, \normHs{1}{\Theta}$, and $T$. 
    In particular, 
    $$\normH{U(t)-u(t)}+\normH{B(t)-b(t)}+\normL{\Theta(t)-\theta(t)}\rightarrow 0\,\,\,\text{as}\,\,\,\epsilon\to 0^{+}.$$ 
\end{theorem}

This paper is organized as follows. In Section~\ref{exist_reg}, we provide all {\it a priori} estimates for the global existence and uniqueness of System~(\ref{Sys1}), as well as the proof of the higher-order regularity of the solution.; while in Section~\ref{conv}, we show the convergence of the solution to System~(\ref{Sys1}) to that of System~(\ref{Sys2}). We also point out that the constant $C$, and others that appear frequently in this paper, may sometimes be an absolute constant that depends only on the domain and the dimension of the relevant problem or quantities, or, may sometimes depend on the parameters in concern; and for the sake of simplicity, we omit all the explicit dependence, and note that they can be traced whenever necessary. 

\section{Proof of Global Existence and Regularity Results}\label{exist_reg}
In this section, we provide the {\it a priori} estimates in order to obtain the global well-posedness of the solution to system \eqref{Sys1}, as well as the higher-order regularity of the solution. Note that the solution can be constructed rigorously by approximations with Galerkin ODEs, following the approach in \cite{Larios_Pei_2017_MHDB_PS} and the reference therein. However, for the sake of simplicity, we omit the details here. 

\subsection{Existence of Solutions to System~(\ref{Sys1}) (Theorem \ref{thm_globExi})}\label{well-posed}
We prove Theorem \ref{thm_globExi} using Galerkin approximation methods. Let $Q_n$ denote the Galerkin projection onto the first $n$ eigenmodes of the Laplacian operator $-\lap$, and let $P_n = P_\sigma Q_n$ denote the projection onto the first $n$ eigenmodes of the Stokes operator $A$. We will use $(u^n, b^n, \theta^n)$ to denote a solution to the Galerkin system, written below in functional form:
\begin{equation}
\left\{ \begin{aligned}
& \frac{d}{dt}u^n + \nu Au^n + P_n\mathcal{B}(u^n, u^n) 
= P_n\mathcal{B}(b^n, b^n) + gP_\sigma \theta^n \vec{e}_2,\\ &
\frac{d}{dt}b^n + \mu Ab^n + P_n\mathcal{B}(u^n, b^n) 
= P_n\mathcal{B}(b^n, u^n), \\ &
\frac{d}{dt}\theta^n -\kappa\lap\theta^n + Q_n \mathcal{B}(u^n, \theta^n) =
\alpha \mu Q_n \lp \bze(|\nabla\times b^n|)|\nabla\times b^n| \rp,
\\&
u^n(0) = P_n(u_0), \quad 
b^n(0) = P_n(b_0), \quad 
\theta^n(0) = Q_n(\theta_0).
\end{aligned}
\right.
\label{Sys_Galerkin}
\end{equation}

In this system, the time evolution equations $u^n$ and $b^n$ are both finite-dimensional ODEs on $P_n(H)$ and the time evolution of $\theta^n$ is a finite-dimensional ODE on $Q_n(L^2(\T))$. Moreover, each time derivative is given by a locally Lipschitz function (see, e.g., \cite{Enlow_Larios_Wu_2023_KSE, Enlow_Larios_Wu_2023_NSE}), hence the short-time existence and uniqueness of the solution $(u^n, b^n, \theta^n)$ is known up to a maximum time interval of existence $[0,T_n)$. To deduce the global existence of solutions for \ref{Sys1}, we work to obtain bounds on $u^n$, $b^n$, and $\theta^n$ that are independent of $n$. In the arguments that follow, we write $u^n \equiv u$, $b^n \equiv b$, and $\theta^n \equiv \theta$ and omit any projection operators for notational simplicity.

\subsubsection{$L^2$-estimates of Theorem~\ref{thm_globExi}} 
\label{thm_main_L2_est}
Multiply the three equations in \ref{Sys_Galerkin} by $u, b, \theta$, respectively, integrate by parts over $\nT^2$, and add, so that we obtain 
\begin{align*}
    &
    \frac{1}{2}\frac{d}{dt}\left(\normL{u}^2+\normL{b}^2+\normL{\theta}^2\right)
    +
    \nu\normL{\nabla u}^2 
    + 
    \mu\normL{\nabla b}^2 
    +
    \kappa\normL{\nabla \theta}^2
    \\&
    =
    \int_{\nT^2}g\theta\vec{e}_2\cdot u\,dx
    +
    \int_{\nT^2}\alpha\mu \bze(|\nabla\times b|)|\nabla\times b|\theta\,dx
    \\&
    \leq
   \frac{g}{2}\normL{u}^2 
   +
   \frac{g}{2}\normL{\theta}^2
   +
   \frac{\mu}{2}\normL{\nabla b}^2
   +
   \frac{\mu\alpha^2 M_{\epsilon}^2}{2}\normL{\theta}^2,
\end{align*}
where we used (\ref{enstrophy_miracle}), (\ref{jacobi}), the divergence-free condition, Young' inequality, and the boundedness of $\bze$. 
Then, integrating in time from $0$ to $T>0$, and by Gr\"onwall inequality, we get for all $t\in [0,T]$,
\begin{align}\label{L2-bound}
    &
    \normL{u(t)}^2
    +
    \normL{b(t)}^2
    +
    \normL{\theta(t)}^2
    \\&\quad
    +
    \nu\int_{0}^{T}\normL{\nabla u(t)}^2\,dt
    +
    \mu\int_{0}^{T}\normL{\nabla b(t)}^2\,dt
    +
    \kappa\int_{0}^{T}\normL{\nabla \theta(t)}^2\,dt \nonumber
    \\&
    \leq
    K_1,  \nonumber
\end{align}
where $K_1$ is a constant that depends on the initial data, as well as $g,\mu, \alpha, \epsilon$,  and $T$,
and does not depend on $n$.

\subsubsection{$H^1$-estimates of Theorem~\ref{thm_globExi}}
\label{thm_main_H1_est}
Multiply the three equations in \ref{Sys_Galerkin} by $A u, Ab, -\Delta\theta$, respectively, integrate by parts over $\nT^2$, and add, so that we obtain 
\begin{align*}
    &
     \frac{1}{2}\frac{d}{dt}\left(\normL{\nabla u}^2+\normL{\nabla b}^2+\normL{\nabla \theta}^2\right)
    +
    \nu\normL{\Delta u}^2 
    + 
    \mu\normL{\Delta b}^2 
    +
    \kappa\normL{\Delta \theta}^2
    \\&
    =
    -
    \int_{\nT^2}(b\cdot\nabla)b\cdot\Delta u\,dx
    +
    \int_{\nT^2}(u\cdot\nabla)b\cdot\Delta b\,dx
    -
    \int_{\nT^2}(b\cdot\nabla)u\cdot\Delta b\,dx
    \\&\quad\quad
    -
    g\int_{\nT^2}\theta \vec{e}_2\cdot\Delta u\,dx
    +
    \int_{\nT^2}(u\cdot\nabla)\theta\Delta\theta\,dx
    \\&\quad\quad\quad
    -
    \alpha\mu\int_{\nT^2}\bze(|\nabla\times b|)|\nabla\times b|\Delta\theta\,dx.
\end{align*}
Next, we estimate the six terms on the right side of the above equations. 
First, by (\ref{jacobi}), the sum of the first three terms is simplified, and thence estimated as 
\begin{align*}
    &
    -
    \int_{\nT^2}(b\cdot\nabla)b\cdot\Delta u\,dx
    +
    \int_{\nT^2}(u\cdot\nabla)b\cdot\Delta b\,dx
    -
    \int_{\nT^2}(b\cdot\nabla)u\cdot\Delta b\,dx
    \\&
    =
    2\int_{\nT^2}(u\cdot\nabla)b\cdot\Delta b\,dx
    \leq
    C\fournorm{u}\fournorm{\nabla b}\normL{\Delta b}
    \\&
    \leq
    C\normL{u}^{1/2}(\normL{u} + \normL{\nabla u})^{1/2}\normL{\nabla b}^{1/2}\normL{\Delta b}^{3/2}
    \\&
    \leq
    C\normL{\nabla b}^2
    +
    C\normL{\nabla u}^2\normL{\nabla b}^2
    +
    \frac{\mu}{4}\normL{\Delta b}^2,
\end{align*}
where we used (\ref{ladyzhenskaya}), Young's inequality, and the $L^2$-bounds obtained in (\ref{L2-bound}). 
Regarding the fourth term, we integrate by parts, apply Young's inequality, and get 
\begin{align*}
    -
    g\int_{\nT^2}\theta \vec{e}_2\cdot\Delta u\,dx
    \leq
    \frac{g}{2}\normL{\nabla\theta}^2
    +
    \frac{g}{2}\normL{\nabla u}^2.
\end{align*}
The fifth term is estimated similarly to that of the first three, and we have
\begin{align*}
    \int_{\nT^2}(u\cdot\nabla)\theta\Delta\theta\,dx
    &
    \leq
    \int_{\nT^2}|u||\nabla\theta||\Delta\theta|\,dx
    \leq
    C\fournorm{u}\fournorm{\nabla\theta}\normL{\Delta\theta}
    \\&
    \leq
    C\normL{u}^{1/2}(\normL{u} + \normL{\nabla u})^{1/2}\normL{\nabla\theta}^{1/2}\normL{\Delta\theta}^{3/2}
    \\&
    \leq
    C\normL{\nabla\theta}^2
    +
    C\normL{\nabla u}^2\normL{\nabla\theta}^2
    +
    \frac{\kappa}{4}\normL{\Delta\theta}^2,
\end{align*}
where we also used the $L^2$-bounds of $u$ obtained in (\ref{L2-bound}). 

As for the last term, by the boundedness of $\bze$, and Young's inequality, we get 
\begin{align*}
    -
    \alpha\mu\int_{\nT^2}\bze(|\nabla\times b|)|\nabla\times b|\Delta\theta\,dx
    &
    \leq
    \frac{\alpha^2\mu^2 M_{\epsilon}^2}{\kappa}\normL{\nabla b}^2
    +
    \frac{\kappa}{4}\normL{\Delta\theta}^2.
\end{align*}
Therefore, by combining all the above estimates, after some rearrangement and simplification, we obtain 
\begin{align*}
     &
     \frac{d}{dt}\left(\normL{\nabla u}^2+\normL{\nabla b}^2+\normL{\nabla \theta}^2\right)
    +
    \nu\normL{\Delta u}^2 
    + 
    \mu\normL{\Delta b}^2 
    +
    \kappa\normL{\Delta \theta}^2
    \\&
    \leq
    \frac{g}{2}\normL{\nabla u}^2 
    +
    \frac{g}{2}\normL{\nabla\theta}^2
    +
    \frac{\alpha^2\mu^2 M_{\epsilon}^2}{\kappa}\normL{\nabla b}^2
    \\&\quad
    +
    C\normL{\nabla u}^2\big( \normL{\nabla b}^2 + \normL{\nabla\theta}^2\big).
\end{align*}
Observing the $L^2$-integrability in time of $\normL{\nabla u}$ from \eqref{L2-bound}, and by the Gr\"onwall inequality, we integrate the above inequality in time from $0$ to $T$, $T>0$, to obtain
\begin{align}\label{H1-bound}
    &
    \normL{\nabla u(T)}^2
    +
    \normL{\nabla b(T)}^2
    +
    \normL{\nabla\theta(T)}^2
    \\&\quad
    +
    \nu\int_{0}^{T}\normL{\Delta u(t)}^2\,dt
    +
    \mu\int_{0}^{T}\normL{\Delta b(t)}^2\,dt
    +
    \kappa\int_{0}^{T}\normL{\Delta\theta(t)}^2\,dt \nonumber
    \\&
    \leq
    K_2,  \nonumber
\end{align} 
where the constant $K_2$ depends on those relevant parameters, as well as on $T$ and $K_1$ in (\ref{L2-bound}), and is independent of $n$.  

\subsubsection{Estimates on time derivatives for Theorem \ref{thm_globExi}}

In order to make valid the convergence arguments used in Section \ref{conv} we require $u,b\in C([0,T]; V)$ and $\theta\in C([0,T]; L^2)$, for which we need 
$\dt u, \dt b \in L^2(0,T; H)$ and $\theta \in L^2(0,T; H^{-1})$. 
We will show this by selecting arbitrary test functions $\phi \in L^2(0,T; H)$ and $\psi \in L^2(0,T; H^1)$ and applying the estimates obtained in Sections \ref{thm_main_L2_est} and \ref{thm_main_H1_est}. Note that by estimate \eqref{H1-bound} and Agmon's inequality \eqref{agmon}, we have $u$ and $b$ bounded above by $K_2$ in $ L^2(0,T; L^\infty) \cap L^\infty(0,T; V)$.

First we take the action of $\dt u$ on $\phi$, from which we obtain 
\begin{align}
\label{dt_est1}
    & \abs{\int_0^T \left\langle \dt u, \phi\right\rangle dt} 
    = \abs{\int_0^T \left\langle \nu\lap u - (u\cdot \grad)u + (b\cdot \grad)b + g\theta \vec{e_2}, \phi\right\rangle dt} 
    \\ &\leq \notag 
    \nu \int_0^T \abs{\lp \lap u, \phi\rp} dt +
    \int_0^T \abs{\lp (u\cdot\grad)u, \phi\rp} dt 
    \\ &\quad + \notag
    \int_0^T \abs{\lp (b\cdot\grad)b, \phi\rp} dt +
    g\int_0^T \abs{\lp \theta \vec{e_2}, \phi\rp} dt 
    \\ &\leq \notag 
    \nu \int_0^T\normLp{2}{\lap u}\normLp{2}{\phi} dt +     
    \int_0^T\normLp{\infty}{u}\normLp{2}{\grad u}\normLp{2}{\phi} dt 
    \\ &\quad + \notag
    \int_0^T\normLp{\infty}{b}\normLp{2}{\grad b}\normLp{2}{\phi} dt +
    g \int_0^T\normLp{2}{\theta}\normLp{2}{\phi} dt \\ &\leq \notag
    C \lp 
    K_2 + 2K_2^2 + K_1
    \rp 
    \int_0^T\normLp{2}{\phi}^2 dt,
\end{align}
where we used H\"older's inequality, Cauchy-Schwarz inequality, and the aforementioned estimates. We then take the action of $\dt b$ on $\phi$ and, \textit{mutatis mutandis}, we obtain 
\begin{align}
    \label{dt_est2}
    \abs{\int_0^T \left\langle \dt b, \phi\right\rangle dt} \leq 
    C \lp 
    \mu K_2 + 2K_2^2 
    \rp 
    \int_0^T\normLp{2}{\phi}^2 dt.
\end{align}
For $\dt \theta$, we take the action on $\psi$ and see that 
\begin{align}
\label{dt_est3}
    &\quad  \abs{\int_0^T \la \dt \theta ,\psi \ra dt} \\ \notag
    &= 
    \abs{\int_0^T \la \kappa \lap \theta - u\cdot \grad \theta 
    + 
    \alpha \mu \bze(\abs{\grad \times b})\abs{\grad \times b} , \psi \ra dt} 
    \\ &\leq \notag 
    \int_0^T \kappa \normLp{2}{\grad\theta}\normLp{2}{\grad\psi} dt + 
    \int_0^T \normLp{\infty}{u}\normLp{2}{\grad \theta}\normLp{2}{\psi} dt \\ & + \notag
    \alpha\mu\int_0^T \normLp{\infty}{\bze}\normLp{2}{\grad b}\normLp{2}{\psi} dt \\ &\leq 
    C \lp 
    K_1 + K_2^2 + \alpha K_1 \normLp{\infty}{\bze}
    \rp
    \int_0^T \normHs{1}{\psi}^2 dt. \notag
\end{align}
By the Aubin-Lions Theorem we deduce that 
\begin{align}
   u, b \in C([0,T]; V) 
\quad\text{ and }\quad
\theta \in C([0,T]; L^2). 
\end{align}

\subsubsection{Convergence of the Galerkin system \eqref{Sys_Galerkin} to the calmed system \eqref{Sys1}}
We now use the previous estimates to justify that our solution to the Galerkin system, now being referred to as $(u^n, b^n, \theta^n)$, converges to a solution of system \eqref{Sys1}. Using the $n$-independent bounds obtained in \eqref{L2-bound}, \eqref{H1-bound}, \eqref{dt_est1}, \eqref{dt_est2}, \eqref{dt_est3} and applying Banach-Alaoglu and Aubin-Lions, we deduce the existence of a subsequence $(u^n, b^n, \theta^n)$, which we will not relabel, and a limit point $(u, b, \theta)$, such that 
\begin{align}
    \label{gal_convList}
    u^n  \rightharpoonup u\,\text{ and }\, b^n & \rightharpoonup b\quad \text{weakly in } 
    L^2(0,T;H^2 \cap V),\\ \notag
    \theta^n & \rightharpoonup \theta \quad \text{weakly in } 
    L^2(0,T;H^1),\\ \notag
    \dt u^n  \rightharpoonup \dt u
    \,\text{ and }\,\dt b^n & \rightharpoonup \dt b \quad \text{weakly in } 
    L^2(0,T; H),\\ \notag
    \dt \theta^n & \rightharpoonup \dt \theta \quad \text{weakly in } 
    L^2(0,T; H^{-1}),\\ \notag
    u^n  \to u\,\text{ and }\, b^n & \to b \quad \text{strongly in } 
    C([0,T]; V),\\ \notag
    \theta^n & \to \theta \quad \text{strongly in }
    C([0,T]; L^2).
 \end{align}
Now we show that, as $n \to \infty$, the Galerkin system \eqref{Sys_Galerkin} converges to system \eqref{Sys1}. Showing the convergence of most of the nonlinear terms is very standard (for a detailed example, see, e.g., \cite{larios2014higher} and the references therein), so brevity we will only show convergence for the calmed Ohmic heating term. First, we write the term as 
\begin{align}
\label{gal_diff}
    &\quad \alpha\mu \bze(\abs{\grad\times b^n})\abs{\grad\times b^n} 
    \\ &= \notag
    \alpha\mu\bze(\abs{\grad\times b^n})
    \lp \abs{\grad\times b^n} - \abs{\grad\times b}  \rp \\ &\quad+ \notag 
    \alpha\mu \lp \bze(\abs{\grad\times b^n}) -  \bze(\abs{\grad\times b}) \rp
    \abs{\grad\times b} \\ &\quad+ \notag
    \alpha\mu\bze(\abs{\grad\times b}) \abs{\grad\times b}.
\end{align}
We now take the action of \eqref{gal_diff} on $\psi\in L^2(0,T; H^1)$ and integrate in time. For the first term, we obtain
\begin{align}
&\quad \int_0^T 
\abs{
\left\langle 
\alpha\mu\bze(\abs{\grad\times b^n} )
\lp \abs{\grad\times b^n} - \abs{\grad\times b}  \rp, 
\psi
\right\rangle }dt \\ &\leq \notag
\int_0^T
\alpha \mu \normLp{\infty}{\bze}
\normLp{2}{\grad b^n - \grad b  }
\normLp{2}{\psi}dt 
\\ &\leq \notag
\alpha\mu \normLp{\infty}{\bze}
\lp \int_0^T \normLp{2}{\grad b^n - \grad b  }^2 dt \rp^{\frac{1}{2}}
\lp \int_0^T \normLp{2}{\psi}^2 dt \rp^{\frac{1}{2}}.
\end{align}
For the second term, similar techniques yield
\begin{align}\label{gal_est1}
&\quad \int_0^T 
\abs{
\left\langle 
\alpha\mu \lp \bze(\abs{\grad\times b^n}) -  \bze(\abs{\grad\times b}) \rp
\abs{\grad\times b},
\psi
\right\rangle }dt \\ &\leq \notag
\alpha\mu 
\int_0^T 
\normLp{2}{\grad b^n - \grad b} 
\normLp{4}{\grad b}
\normLp{4}{\psi}
dt  \\ &\leq \notag
\alpha\mu 
C
\| \grad b^n - \grad b \|_{L^\infty L^2}
\| \grad b \|_{L^\infty L^2}^{\frac{1}{2}} 
\int_0^T 
\normLp{2}{\lap b}^\frac{1}{2}
\normHs{1}{\psi}
dt  \\ &\leq \notag
\alpha\mu 
C
\| \grad b^n - \grad b \|_{L^\infty L^2}
\| \grad b \|_{L^\infty L^2}^{\frac{1}{2}} 
\| \lap b \|^{\frac{1}{2}}_{L^2 L^2}
\| \psi \|_{L^2 H^1}.
\end{align}
We observe that each of these terms converge to $0$ as $n\to\infty$ due to \eqref{gal_convList}, hence $\alpha\mu \bze(\abs{\grad\times b^n})\abs{\grad\times b^n} $ converges to the third term $$\alpha\mu \bze(\abs{\grad\times b})\abs{\grad\times b}.$$ We conclude that solutions to system \eqref{Sys1} exist on $[0,T]$. \qedsymbol

\begin{remark}
    We use the fact that $b^n$ is bounded uniformly in $L^2(0,T;H^2\cap V)$ to obtain estimate \eqref{gal_est1}, which is a consequence of obtaining $H^1$ estimates for solutions with initial data $(u_0, b_0, \theta_0) \in V\times V \times L^2$. The global existence of solutions to \eqref{Sys1} with initial data $(u_0,b_0,\theta_0) \in H\times H\times L^2$ remains an open question. 
\end{remark}

\subsection{Second-order Regularity of Solutions (Theorem \ref{thm_hor})}\label{regularity}
We prove the $H^2$-order regularity of the solution to System~\ref{Sys1} in two steps, similar to \cite{Larios_Pei_Rebholz_2018}. First, we obtain the $H^2$-regularity of $u$ and $b$; then, we use the bounds on the $H^2$-norm of $u$ and $b$ to prove the higher-order regularity of $\theta$. The main reason is to overcome the difficulty of differentiating the absolute value of $\nabla\times b$.  In this section, to simplify notation, we work formally, but the results can be made rigorous by following Galerkin procedure used above.

To begin, we take the action of the relevant equations in system~\ref{Sys1} on $\Delta^2 u$, $\Delta^2 b$, $\Delta\theta$, respectively, integrate by parts over $\nT^2$, and add the results, to obtain
\begin{align*}
    &\,\,
    \frac{1}{2}\frac{d}{dt}\left(\normL{\Delta u}^2+\normL{\Delta b}^2+\normL{\nabla \theta}^2\right)
    \\&\quad
    +
    \nu\normL{\nabla\Delta u}^2 
    + 
    \mu\normL{\nabla\Delta b}^2 
    +
    \kappa\normL{\Delta \theta}^2
    \\&
    =
    g\int_{\nT^2}\theta \vec{e}_2\cdot\Delta^2 u\,dx
    -
    \int_{\nT^2}(u\cdot\nabla)u\cdot\Delta^2 u\,dx
    +
    \int_{\nT^2}(b\cdot\nabla)b\cdot\Delta^2 u\,dx
    \\&\quad
    +
    \int_{\nT^2}(b\cdot\nabla)u\cdot\Delta^2 b\,dx
    -
    \int_{\nT^2}(u\cdot\nabla)b\cdot\Delta^2 b\,dx
    -
    \int_{\nT^2}(u\cdot\nabla)\theta\Delta\theta\,dx
    \\&\quad\quad
    +
    \alpha\mu\int_{\nT^2}\bze(|\nabla\times b|)|\nabla\times b|\Delta\theta\,dx.
\end{align*}
Then, we estimate the seven terms on the right side of the above equation. 
For the first term, integrating by parts twice and applying Young's inequality, we have
\begin{align}\label{hor_est1}
    g\int_{\nT^2}\theta \vec{e}_2\cdot\Delta^2 u\,dx
    \leq
    \frac{g}{2}\normL{\Delta u}^2
    +
    \frac{g}{2}\normL{\Delta \theta}^2. 
\end{align}

The estimates of the remaining five terms are similar, so for the sake of brevity, we provide only the key steps without further clarification. 
Specifically, the second term is estimated as follows,
\begin{align} \label{hor_est2}
    -
    \int_{\nT^2}(u\cdot\nabla)u\cdot\Delta^2 u\,dx
    &
    \leq
    \int_{\nT^2}|\nabla u|^2|\nabla\Delta u|\,dx
    +
    \int_{\nT^2}|u||\nabla\nabla u||\nabla\Delta u|\,dx
    \\& \notag
    \leq
    C\fournorm{\nabla u}^2\normL{\nabla\Delta u}    
    +
    C\fournorm{u}\fournorm{\nabla\nabla u}\normL{\nabla\Delta u}
    \\& \notag
    \leq
    C\normL{\nabla u}\normL{\Delta u}\normL{\nabla\Delta u}
    \\&\quad \notag
    +
    C\normL{u}^{1/2}(\normL{u} + \normL{\nabla u})^{1/2}\normL{\Delta u}^{1/2}\normL{\nabla\Delta u}^{3/2}
    \\&\leq \notag 
     CK_2 \normL{\Delta u}\normL{\nabla\Delta u}
    \\&\quad \notag
    +
    CK_1^{1/2}(K_1 + K_2)^{1/2}\normL{\Delta u}^{1/2}\normL{\nabla\Delta u}^{3/2} 
    \\ &\leq \notag 
    C\lp K_2^2 + K_1^2 \lp K_1 + K_2 \rp^2\rp \normLp{2}{\lap u}^2  + \frac{\nu}{4}\normL{\nabla\Delta u}^2,
\end{align}
where we used the bounds on $\normL{u}$ and $\normL{\nabla u}$ in (\ref{L2-bound}) and (\ref{H1-bound}), as well as (\ref{ladyzhenskaya}) and Young's inequality. 

As for the third term, we proceed somewhat similarly to obtain
\begin{align}\label{hor_est3}
    \int_{\nT^2}(b\cdot\nabla)b\cdot\Delta^2 u\,dx
    &
    \leq
    \int_{\nT^2}|\nabla b|^2|\nabla\Delta u|\,dx
    +
    \int_{\nT^2}|b||\nabla\nabla b||\nabla\Delta u|\,dx
    \\& \notag
    \leq
    \fournorm{\nabla b}^2\normL{\nabla\Delta u}
    +
    \fournorm{b}\fournorm{\nabla\nabla b}\normL{\nabla\Delta u}
    \\& \notag
    \leq
    C\normL{\nabla b}\normL{\Delta b}\normL{\nabla\Delta u}
    \\&\quad \notag
    +
    C\normL{b}^{1/2}(\normL{b} + \normL{\nabla b})^{1/2}\normL{\Delta b}^{1/2}\normL{\nabla\Delta b}^{1/2}\normL{\nabla\Delta u}
    \\& \notag
    \leq 
    C\lp K_2^2 + K_1^2 \lp K_1 + K_2 \rp^2 \rp  \normLp{2}{\lap b}^2
    \\&\quad \notag
    +
    \frac{\nu}{8}\normL{\nabla\Delta u}^2
    +
    \frac{\mu}{8}\normL{\nabla\Delta b}^2,
\end{align}
where we used bounds on both $L^2$- and $H^1$-norms of $b$ obtained in (\ref{L2-bound}) and (\ref{H1-bound}). 
Regarding the fourth term, analogously, by the bounds quoted above, we estimate 
\begin{align}\label{hor_est4}
    \int_{\nT^2}(b\cdot\nabla)u\cdot\Delta^2 b\,dx
    &
    \leq
    \int_{\nT^2}|\nabla b||\nabla u||\nabla\Delta b|\,dx
    +
    \int_{\nT^2}|b||\nabla\nabla u||\nabla\Delta b|\,dx
    \\& \notag
    \leq
    \fournorm{\nabla u}\fournorm{\nabla b}\normL{\nabla\Delta b}
    +
    \fournorm{b}\fournorm{\nabla\nabla u}\normL{\nabla\Delta b}
    \\ &\leq \notag 
    CK_2 \normLp{2}{\lap u}^{1/2}\normLp{2}{\lap b}^{1/2}\normL{\nabla\Delta b} 
    \\ &\quad \notag
    + 
    CK_1^{1/2}\lp K_1 + K_2 \rp^{1/2} \normLp{2}{\lap u}^{1/2}\normL{\nabla\Delta u}^{1/2}\normL{\nabla\Delta b} 
    \\ &\leq \notag 
    C\lp K_2^4 + K_1^2 \lp K_1 + K_2 \rp^2 \rp \normLp{2}{\lap u}^2 + 
    C \normLp{2}{\lap b}^2 
    \\ &\quad \notag +
    \frac{\nu}{8}\normL{\nabla\Delta u}^2
    +
    \frac{\mu}{8}\normL{\nabla\Delta b}^2,
\end{align}
where Young's inequality and (\ref{ladyzhenskaya}) are also used. 

Analogously, for the fifth term, we estimate 
\begin{align}\label{hor_est5}
    -
    \int_{\nT^2}(u\cdot\nabla)b\cdot\Delta^2 b\,dx
    &
    \leq
    C\normLp{2}{\lap u}^2 
    + 
    C\lp K_2^4 + K_1^2 \lp K_1 + K_2\rp^2 \rp \normLp{2}{\lap b}^2 
    \\&\quad \notag
    +
    \frac{\mu}{8} \normLp{2}{\grad\lap b}^2.
\end{align}

Next, we estimate the sixth term as 
\begin{align*}
    -
    \int_{\nT^2}(u\cdot\nabla)\theta\Delta\theta\,dx
    &
    \leq
    C\fournorm{u}\normL{\nabla\theta}^{1/2}\normL{\Delta\theta}^{3/2}
    \\&
    \leq
    CK_1^{1/2} \lp K_1 + K_2 \rp^{1/2} \normLp{2}{\grad \theta}^{1/2} \normLp{2}{\lap \theta}^{3/2} \\ &\leq 
    CK_1^2 \lp K_1 + K_2 \rp^2 \normLp{2}{\grad\theta}^2 + \frac{\kappa}{8} \normLp{2}{\lap\theta}^2,
\end{align*}
where we used the boundedness of $\|u\|_{L^\infty(0,T;H^1)}$ from \eqref{H1-bound}. 

Finally, to estimate the last term, we take advantage of the boundedness of $|\bze|$, apply Poincar\'e's inequality, and obtain 
\begin{align*}
     \alpha\mu\int_{\nT^2}\bze(|\nabla\times b|)|\nabla\times b|\Delta\theta\,dx
     &
     \leq
     \alpha\mu M_{\epsilon}\normL{\nabla\times b}\normL{\Delta\theta}
     \\&
     \leq
     C\alpha^2\mu^2M_{\epsilon}^2\normL{\Delta b}^2 
     +
     \frac{\kappa}{8}\normL{\Delta\theta}^2.
\end{align*}

Now, combining all the above estimates, rearranging and simplifying some terms, we obtain 
\begin{align*}
     &
     \frac{d}{dt}\left(\normL{\Delta u}^2+\normL{\Delta b}^2+\normL{\nabla \theta}^2\right)
    \\&\quad
    +
    \nu\normL{\nabla\Delta u}^2 
    + 
    \mu\normL{\nabla\Delta b}^2 
    +
    \kappa\normL{\Delta \theta}^2
    \\&
    \leq
    \bar{K}\left(\normL{\Delta u}^2+\normL{\Delta b}^2+\normL{\nabla \theta}^2\right),
\end{align*}
where the constant $\bar{K}$ depends on $g,\nu,\mu,\kappa, M_{\epsilon}, \alpha, \lambda_1$, as well as $K_1$ and $K_2$. 
Therefore, integrate the above inequality in time on $[0,T]$, and by Gr\"onwall inequality, we have 
\begin{align}\label{H2-bound-u-b}
    &
    \normL{\Delta u(T)}^2
    +
    \normL{\Delta b(T)}^2
    +
    \normL{\nabla\theta(T)}^2
    \\&\quad
    +
    \nu\int_{0}^{T}\normL{\nabla\Delta u(t)}^2\,dt
    +
    \mu\int_{0}^{T}\normL{\nabla\Delta b(t)}^2\,dt
    +
    \kappa\int_{0}^{T}\normL{\Delta\theta(t)}^2\,dt \nonumber
    \\&
    \leq
    K_3, \nonumber
\end{align} 
where the constant $K_3$ depends on all the aforementioned parameters and on $K_1$ and $K_2$. 

Next, we prove the $H^2$-regularity of $\theta$. Proceeding similarly, we multiply the relevant equations in System~\ref{Sys1} by $\Delta^2 u$, $\Delta^2 b$, $\Delta^2\theta$, respectively, integrate by parts over $\nT^2$, add, and get 
\begin{align*}
    &
    \frac{1}{2}\frac{d}{dt}\left(\normL{\Delta u}^2+\normL{\Delta b}^2+\normL{\Delta \theta}^2\right)
    \\&\quad
    +
    \nu\normL{\nabla\Delta u}^2 
    + 
    \mu\normL{\nabla\Delta b}^2 
    +
    \kappa\normL{\nabla\Delta \theta}^2
    \\&
    =
    g\int_{\nT^2}\theta \vec{e}_2\cdot\Delta^2 u\,dx
    -
    \int_{\nT^2}(u\cdot\nabla)u\cdot\Delta^2 u\,dx
    +
    \int_{\nT^2}(b\cdot\nabla)b\cdot\Delta^2 u\,dx
    \\&\quad
    +
    \int_{\nT^2}(b\cdot\nabla)u\cdot\Delta^2 b\,dx
    -
    \int_{\nT^2}(u\cdot\nabla)b\cdot\Delta^2 b\,dx
    -
    \int_{\nT^2}(u\cdot\nabla)\theta\Delta^2\theta\,dx
    \\&\quad\quad
    +
    \alpha\mu\int_{\nT^2}\bze(|\nabla\times b|)|\nabla\times b|\Delta^2\theta\,dx.
\end{align*}
Then, we estimate the seven terms on the right side of the above equations. 
We will bound the first five terms using the same estimates derived in \eqref{hor_est1}, \eqref{hor_est2}, \eqref{hor_est3}, \eqref{hor_est4}, and \eqref{hor_est5}, respectively.

For the sixth term, 
\begin{align*}
    -
    \int_{\nT^2}(u\cdot\nabla)\theta\Delta^2\theta\,dx
    &
    \leq
    \int_{\nT^2}|\nabla u||\nabla \theta||\nabla\Delta \theta|\,dx
    +
    \int_{\nT^2}|u||\nabla\nabla \theta||\nabla\Delta \theta|\,dx
    \\&
    \leq
    C\fournorm{\nabla u}\fournorm{\nabla \theta}\normL{\nabla\Delta \theta}    
    +
    C\infnorm{u}\normL{\nabla\nabla \theta}\normL{\nabla\Delta \theta}
    \\ &\leq 
    CK_2 \normLp{2}{\lap u}^{1/2}\normLp{2}{\lap \theta}^{1/2} 
    \normLp{2}{\grad\lap \theta} 
    \\&\quad
    +
    C\lp K_1 + K_3 \rp \normLp{2}{\lap\theta} \normLp{2}{\grad\lap \theta} 
    \\ &\leq 
    C\normLp{2}{\lap u}^2 
    + 
    C\lp K_2^4 + \lp K_1 + K_3 \rp^2 \rp \normLp{2}{\lap \theta}^2 
    + 
    \frac{\kappa}{8}\normL{\nabla\Delta \theta}^2,
\end{align*}
where we used the boundedness of $\normL{\nabla u}$ and $\normL{\nabla\theta}$ in \eqref{H1-bound}, and we used \eqref{H2-bound-u-b} along with Agmon's inequality to bound $\normLp{\infty}{u}$ by $C\lp K_1 + K_3\rp$. 
For the last term, we integrate by parts first, then take advantage of the boundedness of $|\bze|$ and $|\nabla\bze|$, and obtain 
\begin{align*}
     &\quad
     \alpha\mu\int_{\nT^2}\bze(|\nabla\times b|)|\nabla\times b|\Delta^2\theta\,dx
     \\&
     \leq
     \alpha\mu M_{\epsilon}\int_{\nT^2}\big|\nabla|\nabla\times b|\big||\nabla\Delta \theta|\,dx
     \\&\quad
     +
     \alpha\mu \int_{\nT^2}|\nabla\bze(|\nabla\times b|)|\big|\nabla|\nabla\times b|\big||\nabla\times b||\nabla\Delta \theta|\,dx
\end{align*}
of which the first term is bounded by 
\begin{align*}
     C\alpha^2 \mu^2 M_{\epsilon}^2\normL{\Delta b}^2
     +
     \frac{\kappa}{8}\normL{\nabla\Delta \theta}^2    
\end{align*}
due to Cauchy-Schwarz inequality and Lemma \ref{evans_lem};
and the second term is bounded by 
\begin{align*}
     &\quad
     C\alpha\mu M_{\epsilon} \fournorm{\nabla b}\fournorm{\nabla\nabla b}\normL{\nabla\Delta\theta}
     \\&
     \leq
     C\alpha \mu M_\epsilon  K_1^{1/2} K_2^{1/2} \normLp{2}{\lap b}^{1/2} \normLp{2}{\grad \lap b}^{1/2} \normLp{2}{\grad \lap \theta} 
     \\&
     \leq
     C\alpha^4\mu^2 M_{\epsilon}^4 K_1^2 K_2^2\normL{\Delta b}^2
     +
     \frac{\mu}{8}\normL{\nabla\Delta b}^2
     +
     \frac{\kappa}{8}\normL{\nabla\Delta \theta}^2.     
\end{align*} 
Now, combining all the above estimates, rearranging and simplifying some terms, we obtain 
\begin{align*}
     &\,\,
     \frac{d}{dt}\left(\normL{\Delta u}^2+\normL{\Delta b}^2+\normL{\Delta \theta}^2\right)
    \\&\quad
    +
    \nu\normL{\nabla\Delta u}^2 
    + 
    \mu\normL{\nabla\Delta b}^2 
    +
    \kappa\normL{\nabla\Delta \theta}^2
    \\&
    \leq
    K^{*}\left(\normL{\Delta u}^2+\normL{\Delta b}^2+\normL{\Delta \theta}^2\right),
\end{align*}
where the constant $K^{*}$ depends on $g,\nu,\mu,\kappa, M_{\epsilon}, \alpha$, as well as $K_1, K_2,$ and $K_3$. 
Finally, integrate the above inequality in time from $0$ to some $T>0$, and by Gr\"onwall inequality, we get 
\begin{align}\label{H2-bound}
    &
    \normL{\Delta u(T)}^2
    +
    \normL{\Delta b(T)}^2
    +
    \normL{\Delta\theta(T)}^2
    \\&\quad
    +
    \nu\int_{0}^{T}\normL{\nabla\Delta u(t)}^2\,dt
    +
    \mu\int_{0}^{T}\normL{\nabla\Delta b(t)}^2\,dt
    +
    \kappa\int_{0}^{T}\normL{\nabla\Delta\theta(t)}^2\,dt \nonumber
    \\&
    \leq
    K_4,  \nonumber
\end{align} 
where the constant $K_4$ depends on all the aforementioned parameters and on $K_1, K_2$, and $K_3$. Thus we have demonstrated the $H^2$-regularity of our solution $(u,b,\theta)$.  \qedsymbol

\begin{remark}
    The constants $K_1, K_2, K_3,$ and $K_4$ all depend on the parameter $M_\epsilon$ stated in \eqref{calm_bounded}, which tends toward infinity as $\epsilon \to 0^+$. Therefore these estimates do not hold as $\epsilon \to 0^+$. 
\end{remark}

A consequence of Theorems \ref{thm_globExi} and \ref{thm_hor} is the global existence of strong solutions on $[0,T]$ which satisfy Definition \ref{def2}. In particular, this result follows from estimate \eqref{H2-bound}.

\subsection{Uniqueness of Strong Solutions (Theorem \ref{thm_uni})}
Following the well-known weak-strong uniqueness argument for the Navier-Stokes equations, it suffices to show that the strong solution is unique. Note that due to the Ohmic heating effect in System~\ref{Sys1}, one needs to work with at least $H^2$ initial condition and regularity of the solution. To begin, we assume there are two distinct strong solutions $(u_{m}, b_{m}, \theta_{m}, p_{m}), m=1,2$, to System~\ref{Sys1}, in the sense of Theorem~\ref{thm_globExi}, with initial data $u_{m,0}$, $b_{m,0}$, and $\theta_{m,0}$. 
Subtract the corresponding equations satisfied by the two solutions, and denote the differences by $\tilde{u} = u_1-u_2$, $\tilde{b}=b_1-b_2$, $\tilde{\theta}=\theta_1-\theta_2$, and $\tilde{p}=p_1-p_2$, as well as $\tilde{u}_0 = u_{1,0}-u_{2,0}$, $\tilde{b}_0=b_{1,0}-b_{2,0}$, $\tilde{\theta}_0=\theta_{1,0}-\theta_{2,0}$. We obtain the system for $\tu$, $\tb$, $\ttheta$, and $\tilde{p}$ as 
\begin{equation}
       \left\{
       \begin{aligned}
       &\frac{\partial \tu}{\partial t} 
       - 
       \nu\Delta \tu 
       + 
       (u_1\cdot\nabla) \tu 
       +
       (\tu\cdot\nabla) u_2 
       +
       \nabla \tilde{p} 
       \\&\quad\quad\quad
       = 
       (b_1\cdot\nabla) \tb
       +
       (\tb\cdot\nabla) b_2
       +
       g\ttheta \vec{e}_2,
       \\&
       \frac{\partial \tb}{\partial t} 
       - 
       \mu\Delta \tb 
       + 
       (u_1\cdot\nabla) \tb
       +
       (\tu\cdot\nabla) b_2
       = 
       (b_1\cdot\nabla) \tu
       +
       (\tb\cdot\nabla) u_2,
       \\&
       \frac{\partial \ttheta}{\partial t} 
       -
        \kappa\Delta\ttheta
       + 
       (u_1\cdot\nabla) \ttheta
       +
       (\tu\cdot\nabla) \theta_2  
       \\&\quad\quad\quad
       = 
       \alpha\mu\big( \bze(|\nabla\times b_1|)|\nabla\times b_1| -  \bze(|\nabla\times b_2|)|\nabla\times b_2| \big), 
    \\&
    \nabla \cdot \tu = 0 = \nabla \cdot \tb,
    \\&
    \tu(0)=\tu_0,\,\, \tb(0) = \tb_0, \,\,\ttheta(0) = \tilde{\theta}_0. 
\end{aligned}
    \right.
    \label{Sys3}
  \end{equation}
Then, multiply $\tu, \tb, \ttheta$ to the relevant equations in system~\ref{Sys3}, respectively, integrate by parts over $\nT^2$, and add, so that we obtain
\begin{align*}
    &\quad
     \frac{1}{2}\frac{d}{dt}\left(\normL{\tu}^2+\normL{\tb}^2+\normL{\ttheta}^2\right)
    +
    \nu\normL{\nabla\tu}^2 
    + 
    \mu\normL{\nabla\tb}^2 
    +
    \kappa\normL{\nabla\ttheta}^2
    \\&
    =
    g\int_{\nT^2}\ttheta \vec{e}_2\cdot \tu\,dx
    -
    \int_{\nT^2}(\tu\cdot\nabla)u_2\cdot \tu\,dx
    +
    \int_{\nT^2}(\tb\cdot\nabla)b_2\cdot \tu\,dx
    \\&\quad
    -
    \int_{\nT^2}(\tu\cdot\nabla)b_2\cdot \tb\,dx
    +
    \int_{\nT^2}(\tb\cdot\nabla)u_2\cdot \tb\,dx
    -
    \int_{\nT^2}(\tu\cdot\nabla)\theta_2 \ttheta\,dx
    \\&\quad\quad
    +
    \alpha\mu\int_{\nT^2}\big( \bze(|\nabla\times b_1|)|\nabla\times b_1| -  \bze(|\nabla\times b_2|)|\nabla\times b_2| \big)\ttheta\,dx,
\end{align*}
where we used the divergence-free condition $\nabla\cdot\tu = 0 = \nabla\cdot\tb$. 

Next, we estimate the seven terms on the right side of the above equation. 
By Cauchy-Schwarz inequality, the first term is estimated as
\begin{align*}
    g\int_{\nT^2}\ttheta \vec{e}_2\cdot \tu\,dx
    \leq
    \frac{g}{2}\normL{\tu}^2 
    +
    \frac{g}{2}\normL{\ttheta}^2.  
\end{align*}

For the second term, by the $H^1$-boundedness
obtained in (\ref{H1-bound}), we have 
\begin{align*}
    -
    \int_{\nT^2}(\tu\cdot\nabla)u_2\cdot \tu\,dx
    &
    \leq
    C\normL{\nabla u_2}\fournorm{\tu}^2
    \leq
    CK_2\normL{\tu}(\normL{\tu} + \normL{\nabla\tu})
    \\&
    \leq
    CK_2^2\normL{\tu}^2
    +
    \frac{\nu}{8}\normL{\nabla\tu}^2.
\end{align*}
The third through the sixth terms are estimated similarly.  For example, the third term can be estimated by
\begin{align*}
    \int_{\nT^2}(\tb\cdot\nabla)b_2\cdot \tu\,dx
    &
    \leq
    C\normL{\nabla b_2}\fournorm{\tb}\fournorm{\tu}
    \\&
    \leq
    CK_2\normL{\tb}^{1/2}(\normL{\tb} + \normL{\nabla\tb})^{1/2}\normL{\tu}^{1/2}(\normL{\tu} + \normL{\nabla\tu})^{1/2}
    \\&
    \leq
    CK_2\normL{\tb}(\normL{\tb} + \normL{\nabla\tb})
    +
    C\normL{\tu}(\normL{\tu} + \normL{\nabla\tu})
    \\&
    \leq
    C\normL{\tu}^2
    +
    CK_2^2\normL{\tb}^2
    +
    \frac{\nu}{8}\normL{\nabla\tu}^2
    +
    \frac{\mu}{8}\normL{\nabla\tb}^2
\end{align*}
where we used both (\ref{L2-bound}) and (\ref{H1-bound}); and for the fourth one, we also have
\begin{align*}
    -
    \int_{\nT^2}(\tu\cdot\nabla)b_2\cdot \tb\,dx
    &
    \leq
    CK_2^2\normL{\tu}^2
    +
    C\normL{\tb}^2
    +
    \frac{\nu}{8}\normL{\nabla\tu}^2
    +
    \frac{\mu}{8}\normL{\nabla\tb}^2;
\end{align*}
while for the fifth term, we take advantage of the $H^1$-bounds of $u_2$ in (\ref{H1-bound}), and obtain an upper bound as 
\begin{align*}
    \int_{\nT^2}(\tb\cdot\nabla)u_2\cdot \tb\,dx
    &
    \leq
    CK_2^2\normL{\tb}^2
    +
    \frac{\mu}{8}\normL{\nabla\tb}^2;
\end{align*}
as for the sixth term, we proceed similarly and get 
\begin{align*}
    &
    -
    \int_{\nT^2}(\tu\cdot\nabla)\theta_2 \ttheta\,dx
    \leq
    CK_2^2\normL{\tu}^2
    +
    C\normL{\ttheta}^2
    +
    \frac{\nu}{8}\normL{\nabla\tu}^2
    +
    \frac{\kappa}{4}\normL{\nabla\ttheta}^2.
\end{align*}
Now, it remains to estimate the last term that involves the Ohmic heating. 
We start by bounding the magnitude of the first factor inside the integrand. 
Thanks to the reverse triangle inequality and the global Lipschitz condition on $\bze$, we obtain
\begin{align*}
    &
    \big| \bze(|\nabla\times b_1|)|\nabla\times b_1| -  \bze(|\nabla\times b_2|)|\nabla\times b_2| \big|
    \\&
    \leq
     \big| \bze(|\nabla\times b_1|)|\nabla\times b_1| - \bze(|\nabla\times b_1|)|\nabla\times b_2| \big|
     \\&\quad\quad
     +
     \big| \bze(|\nabla\times b_1|)|\nabla\times b_2| - \bze(|\nabla\times b_2|)|\nabla\times b_2| \big|
     \\&
     =
     \big| \bze(|\nabla\times b_1|)|\nabla\times \tb|\big|
     +
     \big| \bze(|\nabla\times b_1|) - \bze(|\nabla\times b_2|) \big||\nabla\times b_2|
     \\&
     \leq
     M_{\epsilon}|\nabla\times \tb|
     +
     L_{\epsilon}|\nabla\times \tb||\nabla\times b_2|.
\end{align*}
Thus, inserting the above estimates to the last term, we obtain
\begin{align*}
    &\quad
    \alpha\mu\int_{\nT^2}\big( \bze(|\nabla\times b_1|)|\nabla\times b_1| -  \bze(|\nabla\times b_2|)|\nabla\times b_2| \big)\ttheta\,dx,
    \\&
    \leq
    \alpha\mu M_{\epsilon}\int_{\nT^2}|\nabla\times\tb||\ttheta|\,dx
    +
    \alpha\mu L_{\epsilon}\int_{\nT^2}|\nabla\times b_2||\nabla\times \tb||\ttheta|\,dx
    \\&
    \leq
    C\normL{\ttheta}^2 
    +
    \frac{\mu}{16}\normL{\nabla\tb}^2
    +
    C\fournorm{\nabla b_2}\normL{\nabla\tb}\fournorm{\ttheta}
    \\&
    \leq
    C\normL{\ttheta}^2 
    +
    \frac{\mu}{16}\normL{\nabla\tb}^2
    +
    C\normL{\nabla \tb}\normL{\ttheta}^{1/2}\normL{\nabla\ttheta}^{1/2}
    \\&
    \leq
    C\normL{\ttheta}^2
    +
    \frac{\mu}{8}\normL{\nabla\tb}^2
    +
    \frac{\kappa}{4}\normL{\nabla\ttheta}^2,
\end{align*}
where in the penultimate step we used the $H^2$-bounds obtained in (\ref{H2-bound}). 

Therefore, by collecting all the above estimates, and after some simplification and rearrangement, we finally obtain 
\begin{align*}
     &\quad
     \frac{d}{dt}\left(\normL{\tu}^2+\normL{\tb}^2+\normL{\ttheta}^2\right)
    +
    \nu\normL{\nabla \tu}^2 
    + 
    \mu\normL{\nabla \tb}^2 
    +
    \kappa\normL{\nabla \ttheta}^2
    \\&
    \leq
    \tilde{K}\big( \normL{\tu}^2 + \normL{\tb}^2 + \normL{\ttheta}^2\big),
\end{align*}
where $\tilde{K}$ depends on $g, \nu, \mu, \kappa, \alpha, M_{\epsilon}, L_{\epsilon}$, as well as $K_m, m=1,2,3,4$.  The
Gr\"onwall inequality then yields continuous dependence on initial data.  Moreover, if the initial data are equal, so that $\tu(0)=\tb(0)=0$ and $\ttheta(0)=0$, it follows that $u_1=u_2$, $b_1=b_2$, and $\theta_1=\theta_2$. Hence, the global well-posedness of system \eqref{Sys1}, is now proved. \qedsymbol

\section{Proof of the Convergence Theorem~\ref{thm_conv}}
\label{conv}
In this section, we show that on the common time-interval of existence, the difference between the solution of System~\ref{Sys1} and that of System~\ref{Sys2}, is on the order of $\epsilon$; and in particular, such error approaches $0$ as $\epsilon$ goes to $0^{+}$. One of the difficulties we face in proving convergence of the calmed system is obtaining bounds on the difference equation that are independent of $\epsilon$. Each upper bound $K_1, K_2, K_3,$ and $K_4$ (obtained in \eqref{L2-bound}, \eqref{H1-bound}, \eqref{H2-bound-u-b}, and \eqref{H2-bound} respectively) has a dependence on $\epsilon$, as each $K_i$ has a factor of the form $e^{CM_\epsilon^2 T}$ (or $e^{CM_\epsilon^4 T}$ in the case of $K_4$), and by assumption, $\lim_{\epsilon \to 0^+} M_\epsilon = \infty$. By applying Lemma \ref{boot} we are able to derive bounds independent of $\epsilon$ and complete the proof of convergence.

\subsection{Proof of Theorem~\ref{thm_conv}}
Set $\tu = U - u$, $\tb = B - b$, $\ttheta = \Theta - \theta$, and $\tilde{p} = P - p$. 
We take the difference of equations \eqref{Sys2} and \eqref{Sys1} to obtain the system
\begin{equation}
        \left\{
        \begin{aligned}
        &\frac{\partial \tu}{\partial t} - \nu \lap \tu + \grad \tilde{p} =         
        \lp \tu \cdot \grad \rp \tu - \lp \tu \cdot \grad \rp U - \lp U \cdot \grad \rp \tu 
        \\&\quad\quad\quad\quad\quad\quad\quad\quad\quad +
        \lp B \cdot \grad \rp \tb - \lp \tb \cdot \grad \rp \tb + \lp \tb \cdot \grad \rp B +
        g\ttheta \vec{e}_2, 
        \\&
        \frac{\partial \tb}{\partial t} - 
        \mu\Delta \tb = 
        \lp B \cdot \grad \rp \tu - \lp \tb \cdot \grad \rp \tu + \lp \tb \cdot \grad \rp U 
        \\&\quad\quad\quad\quad\quad\quad + 
        \lp \tu \cdot \grad \rp \tb - \lp \tu \cdot \grad \rp B - \lp U \cdot \grad \rp \tb,
        \\&
        \frac{\partial \ttheta}{\partial t} -
        \kappa\lap \ttheta = 
        \lp \tu \cdot \grad \rp \ttheta - \lp \tu \cdot \grad \rp \Theta - \lp U \cdot \grad \rp \ttheta 
        \\&\quad\quad\quad\quad\quad\quad + 
        \alpha\mu \lp \jj - \bze\lp \jj \rp \rp \jj
        \\&\quad\quad\quad\quad\quad\quad -
        \alpha\mu\lp \bze\lp \jj\rp - \bze\lp\ej\rp \rp \lp \tj \rp
        \\&\quad\quad\quad\quad\quad\quad +
        \alpha\mu \lp \tj \rp  \bze\lp\jj\rp
        \\&\quad\quad\quad\quad\quad\quad +
        \alpha\mu\lp \bze\lp \jj \rp - \bze\lp \ej \rp \rp \jj      
        \\&
        \nabla \cdot \tu = 0 = \nabla \cdot \tb,
        \\&
        \tu(0)=\tb(0) = 0, \,\,\ttheta(0) = 0. 
\end{aligned}
    \right.
    \label{Sys_convergence}
\end{equation}
Set $T_0 = \frac{1}{2}\lp T+T^* \rp$, for some $T^*$ such that $T<T^*$ and for which the solution $(U,B,\Theta)$ exists on $[0,T^*]$.
First we will remark that since $(U, B, \Theta)$ satisfy the conditions of Assumption \ref{thm_short_original}, there exists $M>0$ such that for a.e. $t\in[0,T_0]$, 
\begin{align}\label{Ohm_original_H2_bound}
    \normHs{2}{U(t)} + \normHs{2}{B(t)} + \normHs{1}{\Theta(t)} \leq M,
\end{align}
hence we can bound these terms above by the constant $M$ when needed. 

Next, we select $\epsilon_0 > 0$ to be sufficiently small so that 
\begin{align}\label{eps_choice}
    CM^{2+2\beta}e^{C\lp M^\frac{4}{3} + M^2 \rp T_0} \epsilon_0^{2\gamma} < \frac{1}{2}
\end{align}

Now we apply Lemma \ref{boot}: let $H(T)$ be the statement
\begin{align}\label{boot_hypothesis}
    \sup_{0<\epsilon<\epsilon_0}
    \sup_{0\leq t \leq T}\lp 
    \normHs{1}{\tu (t)}^2 + 
    \normHs{1}{\tb (t)}^2 + 
    \normLp{2}{\ttheta (t)}^2 \rp < 1,
\end{align}
and let $C(T)$ be the statement  
\begin{align}\label{boot_conclusion}
    \sup_{0<\epsilon<\epsilon_0}
    \sup_{0 \leq t \leq T}\lp 
    \normHs{1}{\tu (t)}^2 + 
    \normHs{1}{\tb (t)}^2 + 
    \normLp{2}{\ttheta (t)}^2 \rp < \frac{1}{2}.
\end{align}
We can infer the following about $H(T)$ and $C(T)$ which will simplify the proof of Lemma \ref{boot}:
\begin{itemize}
    \item From the definition of semi-strong solutions, the hypotheses of Assumption \ref{thm_short_original}, and given that $\tu(0) = \tb(0) = 0$ and $\ttheta(0)=0$, we deduce that \eqref{boot_hypothesis} is valid for some $T >0$. Therefore Condition \ref{lem_d} is immediately satisfied.
    
    \item If $C(T)$ is true, then $H(T)$ is true. Therefore, proving Lemma \ref{boot} is tantamount to proving $H(T)$ is true. 
    \item If $H(s_0)$ and $C(t_0)$ are true, then $H(s)$ and $C(t)$ are true for all $s\in[0,s_0]$ and all $t\in[0,t_0]$. Therefore, we may assume that $t_m$ is a monotone increasing sequence in Statement \ref{lem_c}.

\end{itemize}

First we show that Statement \ref{lem_a} is true. 
Assume that $H(T)$ is valid at $T$ and fix $\epsilon \in (0, \epsilon_0)$. To begin, we take the inner product of the first two equations in system \eqref{Sys_convergence} with $\tu$ and $\tb$ respectively, integrate by parts, and remove any terms that are eliminated using identities \eqref{symm1} and \eqref{symm2}:

\begin{align}\label{Ohm_conv1}
    &\quad \frac{1}{2} \frac{d}{dt}\lp \normLp{2}{\tu}^2 
    + 
    \normL{\tb}^2 \rp 
    + 
    \nu \normL{\grad\tu}^2 
    + 
    \mu \normL{\grad\tb}^2 \\ 
    &= \notag 
    - \lp \lp \tu \cdot \grad\rp U, \tu \rp 
    + 
    \lp \lp \tb \cdot \grad\rp B, \tu \rp 
    +
    \lp \lp \tb \cdot \grad\rp U, \tb \rp 
    \\&\quad\notag
    -\lp \lp \tu \cdot \grad\rp B, \tb \rp 
    + 
    \lp g\ttheta \vec{e}_2, \tu \rp. 
\end{align}

We now take the inner product
of the system with the triplet  $\lp -\lap\tu, -\lap \tb, \ttheta \rp$, integrate by parts, invoke \eqref{enstrophy_miracle} and \eqref{jacobi} where appropriate, and obtain
\begin{align} 
    &\quad \label{Ohm_conv2}
    \frac{1}{2}\frac{d}{dt}
    \left(\normL{\grad\tu}^2+\normL{\grad\tb}^2+\normL{\ttheta}^2\right) +
    \nu\normL{\lap\tu}^2 + 
    \mu\normL{\lap\tb}^2 +
    \kappa\normL{\grad\ttheta}^2 
    \\& \notag
    =
    \lp \lp \tu \cdot \grad \rp   U, \lap \tu \rp +
    \lp \lp   U \cdot \grad \rp \tu, \lap \tu \rp -
    \lp \lp   B \cdot \grad \rp \tb, \lap \tu \rp -
    \lp \lp \tb \cdot \grad \rp   B, \lap \tu \rp 
    \\ &\quad - \notag
    \lp g \ttheta \vec{e_2}, \lap \tu \rp -
    \lp \lp   B \cdot \grad \rp \tu, \lap \tb \rp -
    \lp \lp \tb \cdot \grad \rp   U, \lap \tb \rp -
    2 \lp \lp \tu \cdot \grad \rp \tb, \lap \tb \rp 
    \\ &\quad+ \notag
    \lp \lp \tu \cdot \grad \rp   B, \lap \tb \rp +
    \lp \lp   U \cdot \grad \rp \tb, \lap \tb \rp +
    \lp \lp \tu \cdot \grad \rp \Theta, \ttheta \rp +
    N_1 + N_2 + N_3 + N_4,
\end{align}
where 
\begin{align}
    \label{nonlinear_terms} 
    & N_1 = \alpha\mu \int_{\T} \abs{\jj - \bze\lp \jj \rp}\jj \ttheta dx, \\ \notag
    & N_2 = \alpha\mu \int_{\T} \lp \bze\lp\jj\rp - \bze\lp\ej\rp \rp \lp \tj \rp \ttheta dx, \\ \notag
    & N_3 = \alpha\mu \int_{\T}\lp \tj \rp \bze(\jj) \ttheta dx, \\ \notag 
    & N_4 = \alpha\mu \int_{\T} \lp \bze(\jj) - \bze(\ej) \rp \jj \ttheta dx.
\end{align}
We now combine \eqref{Ohm_conv1} and \eqref{Ohm_conv2} to obtain the equation
\begin{align}
    \label{Ohm_conv3}
    &\quad \frac{1}{2}\frac{d}{dt}
    \left(\normHs{1}{\tu}^2+\normHs{1}{\tb}^2+\normL{\ttheta}^2\right) 
    +
    \nu\normHs{1}{\grad\tu}^2 
    + 
    \mu\normHs{1}{\grad\tb}^2 
    +
    \kappa\normL{\grad\ttheta}^2 
    \\& \notag
    =
    - \lp \lp \tu \cdot \grad\rp U, \tu \rp 
    + 
    \lp \lp \tb \cdot \grad\rp B, \tu \rp 
    +
    \lp \lp \tb \cdot \grad\rp U, \tb \rp 
    -
    \lp \lp \tu \cdot \grad\rp B, \tb \rp 
    \\&\quad\notag 
    +
    \lp g\ttheta \vec{e}_2, \tu \rp 
    +
    \lp \lp \tu \cdot \grad \rp   U, \lap \tu \rp 
    +
    \lp \lp   U \cdot \grad \rp \tu, \lap \tu \rp 
    -
    \lp \lp   B \cdot \grad \rp \tb, \lap \tu \rp 
    \\&\quad\notag
    -
    \lp \lp \tb \cdot \grad \rp   B, \lap \tu \rp 
    -
    \lp g \ttheta \vec{e_2}, \lap \tu \rp 
    -
    \lp \lp   B \cdot \grad \rp \tu, \lap \tb \rp 
    -
    \lp \lp \tb \cdot \grad \rp   U, \lap \tb \rp 
    \\&\quad\notag
    -
    2 \lp \lp \tu \cdot \grad \rp \tb, \lap \tb \rp 
    +
    \lp \lp \tu \cdot \grad \rp   B, \lap \tb \rp 
    +
    \lp \lp   U \cdot \grad \rp \tb, \lap \tb \rp 
    +
    \lp \lp \tu \cdot \grad \rp \Theta, \ttheta \rp 
    \\&\quad\notag 
    +
    N_1 + N_2 + N_3 + N_4.
\end{align}

The first twelve terms can be bounded above in a similar fashion using H\"older's inequality, Agmon's inequality, Young's inequality, and  \eqref{Ohm_original_H2_bound}. For the sake of brevity, we show this only for the first term: 
\begin{align} \label{Ohm_conv4}
    \abs{\lp \lp \tu \cdot \grad \rp U, \tu \rp} &\leq 
    \normLp{\infty}{\tu} \normLp{2}{\grad U} \normLp{2}{\tu} \\ &\leq \notag
    C \normLp{2}{\tu}^{1/2} \normHs{2}{\tu}^{1/2} 
    \normLp{2}{\grad U} \normLp{2}{\tu} \\ &\leq \notag
    C \normLp{2}{\grad U}^{4/3} \normLp{2}{\tu}^2 + \frac{\nu}{22}\normHs{2}{\tu}^2  \\ &\leq \notag
    C M^\frac{4}{3}\normLp{2}{\tu}^2 + \frac{\nu}{22}\normHs{2}{\tu}^2.
\end{align}
The next four terms can be bounded using the aforementioned inequalities and \eqref{boot_hypothesis}. Again we show this only for the first term:
\begin{align} \label{Ohm_conv5}
    2\abs{ \lp \lp \tu \cdot \grad \rp \tb, \lap \tb \rp } &\leq 
    2 \normLp{\infty}{\tu}
    \normLp{2}{\grad\tb}
    \normLp{2}{\lap\tb} \\ &\leq \notag 
    C \normLp{2}{\tu}^{1/2}\normHs{2}{\tu}^{1/2}
    \normLp{2}{\grad\tb}
    \normLp{2}{\lap\tb} \\ &\leq \notag 
    C \normLp{2}{\tu}^{1/2}\normHs{2}{\tu}^{1/2}
    \normLp{2}{\lap\tb} \\ &\leq \notag 
    C\normLp{2}{\tu}^2 + \frac{\nu}{22}\normHs{2}{\tu}^2 + \frac{\mu}{18}\normLp{2}{\lap\tb}^2,
\end{align}
noting that \eqref{boot_hypothesis} was used to obtain the bound $\normLp{2}{\grad\tb} < 1$.

We now focus on upper bounds for each $N_i$: For $N_1$, using inequalities \eqref{calm_convergence_rate}, \eqref{PT1}, and \eqref{Ohm_original_H2_bound}, we obtain

\begin{align}    
    \label{Ohm_conv_n1}
    \abs{N_1} &= \alpha \mu \abs{\int_{\T} \abs{\jj - \bze\lp \jj \rp}\jj \ttheta dx} \\ &\leq \notag
    C \epsilon^\gamma \int_{\T} \abs{\grad B}^{\beta + 1} \ttheta dx \\ &\leq \notag
    C \epsilon^\gamma \normLp{2\beta + 2}{\grad B}^{\beta + 1} \normL{\ttheta} \\ &\leq \notag
    C \epsilon^\gamma \normL{\grad B} \normL{\lap B}^\beta \normL{\ttheta} \\ &\leq \notag
    C \epsilon^{2\gamma}  \normL{\grad B}^2 \normL{\lap B}^{2\beta} + 
    \normL{\ttheta}^2 \\ &\leq \notag
    C M^{2+2\beta} \epsilon^{2\gamma} + 
    \normL{\ttheta}^2.
\end{align}
For $N_2$, using the Lipschitz property of $\bze$, the reverse triangle inequality, \eqref{boot_hypothesis}, and \eqref{ladyzhenskaya}, we obtain
\begin{align}
    \label{Ohm_conv_n2}
   \abs{N_2} &=\alpha\mu \abs{ \int_{\T} \lp \bze\lp\jj\rp - \bze\lp\ej\rp \rp \lp \tj \rp \ttheta dx} \\ \notag
   & \leq C \int_{\T} \abs{\grad \tb}^2 \abs{\ttheta} dx \\ \notag
   & \leq C ||\grad \tb ||_{L^4}^2 || \ttheta ||_{L^2} \\ \notag
   &\leq C||\grad \tb ||_{L^2} 
   ||\lap \tb ||_{L^2}  || \ttheta ||_{L^2} \\ \notag
   &\leq C
   ||\lap \tb ||_{L^2}  || \ttheta ||_{L^2} \\ \notag
  &\leq C || \ttheta ||_{L^2}^2 +
  \frac{\mu}{18}||\lap \tb ||_{L^2}^2,
\end{align}
in which the estimate $\normLp{2}{\grad \tb} < 1$ is obtained from \eqref{boot_hypothesis}.  
For $N_3$, we use \eqref{Ohm_original_H2_bound} and \eqref{boot_hypothesis} which yield
\begin{align} 
    \label{Ohm_conv_n3}
     \abs{N_3} &= \abs{\int_{\T} \alpha \mu \lp \tj \rp \bze(\jj) \ttheta dx} \\ \notag &\leq 
    C \normLp{4}{\grad B} || \grad \tb ||_{L^2} || \ttheta ||_{L^4} \\ \notag &\leq 
    C \normLp{2}{\grad B}^\frac{1}{2}\normLp{2}{\lap B}^\frac{1}{2}
    \normLp{2}{\grad \tb}
    \normLp{2}{\ttheta}^\frac{1}{2}
    \normHs{1}{\ttheta}^\frac{1}{2} \\ \notag &\leq
    C M \normLp{2}{\grad \tb}
    \lp  
    \normLp{2}{\ttheta} + \normLp{2}{\ttheta}^\frac{1}{2}\normLp{2}{\grad \ttheta}^\frac{1}{2} \rp 
    \\ \notag &\leq
    CM^2 \normLp{2}{\grad \tb}^2 + C\normLp{2}{\ttheta}^2 + \frac{\kappa}{4} ||\grad \ttheta ||_{L^2}^2,
\end{align}
and $N_4$ can be bounded in the exact same manner as $N_3$ thanks to the Lipschitz condition on $\bze$. 

Finally, we can bound the left hand side of \eqref{Ohm_conv3} using the estimates in \eqref{Ohm_conv4}, \eqref{Ohm_conv5}, \eqref{Ohm_conv_n1}, \eqref{Ohm_conv_n2}, \eqref{Ohm_conv_n3}, then rearrange the terms to obtain
\begin{align} \label{Ohm_conv6}
&\quad \frac{d}{dt}
\left(\normHs{1}{\tu}^2+\normHs{1}{\tb}^2+\normL{\ttheta}^2\right) +
\nu\normHs{1}{\grad\tu}^2 + 
\mu\normHs{1}{\grad\tb}^2 +
\kappa\normL{\grad\ttheta}^2  \\ &\leq \notag 
C\lp M^\frac{4}{3} + M^2 \rp\left(\normHs{1}{\tu}^2+\normHs{1}{\tb}^2+\normL{\ttheta}^2\right) + 
C M^{2+2\beta} \epsilon^{2\gamma}.
\end{align}
Using Gr\"onwall's inequality \eqref{Gronwall} and the fact that $\tu(0) = \tb(0) = 0$, $\ttheta(0) = 0$, from \eqref{Ohm_conv6} we deduce that for all $t\in [0,T]$,
\begin{align}\label{Ohm_conv_rate1}
    \normHs{1}{\tu(t)}^2+\normHs{1}{\tb(t)}^2+\normL{\ttheta(t)}^2 \leq
    C M^{2+2\beta}e^{C\lp M^{\frac{4}{3}} + M^2 \rp t}
    \epsilon^{2\gamma}.
\end{align}

This estimate establishes the desired convergence. We continue in showing that each condition of Lemma \ref{boot} is satisfied. Take the supremum of both sides of \eqref{Ohm_conv_rate1} over all $\epsilon \in (0, \epsilon_0)$ and all $t \in [0,T]$ and apply \eqref{eps_choice}, so we have
\begin{align}\label{Ohm_conv_rate2}
& \sup_{0<\epsilon<\epsilon_0}
\sup_{0<t<T}\lp \normHs{1}{\tu(t)}^2+\normHs{1}{\tb(t)}^2+\normL{\ttheta(t)}^2 \rp 
\\ &\leq \notag
C M^{2+2\beta}e^{C\lp M^{\frac{4}{3}} + M^2 \rp T_0}
\epsilon_0^{2\gamma} \\ &< \notag
\frac{1}{2}.
\end{align}

Thus we have shown that $C(T)$ is valid, proving Statement \ref{lem_a} of Lemma \ref{boot}. 

To prove Statement \ref{lem_b} we first fix $\delta > 0$ such that $\delta < \frac{1}{2}\lp T_0 - T \rp$. Now suppose that $C(t_0)$ holds for some $t_0 \in [0,T]$. By our choice of $\delta$, we have $t_0 + \delta < T_0$ independent of the value of $t_0$. 
Since the bound stated in Inequality \eqref{Ohm_original_H2_bound} is valid a.e. on $[0, t_0 + \delta]$, it follows that \eqref{Ohm_conv_rate1} and \eqref{Ohm_conv_rate2} are valid for $T \equiv t_0 + \delta$. Therefore $C(t_0 + \delta)$, and thus $H(t_0 + \delta)$, are true, proving Statement \ref{lem_b}.

Similarly, if $C(t_m)$ holds for each $m$ and $t_m$ converges to $t_c$ for some $t_c \in [0,T]$, then there exists $m_0$ such that $t_c \in (t_m, t_m + \delta)$ for all $m \geq m_0$. Thus $C(t_c)$ is true, proving Statement \ref{lem_c}.

Finally, we have already discussed that Statement \ref{lem_d} of Lemma \ref{boot} is true. It follows that $C(t)$ (and thus $H(t)$) is true for all $t\in [0,T_0]$. Therefore the bootstrapping argument has been justified. 

Now, with all the conditions of Lemma \ref{boot} verified, we may apply Inequality \eqref{boot_hypothesis} to derive the resulting inequalities, which yield Inequality \eqref{Ohm_conv_rate1}. From \eqref{Ohm_conv_rate1} we ascertain that Estimate \eqref{thm_conv_rate1} holds, with 
\begin{align}\label{c_1}
    c_1 = C M^{2+2\beta}e^{C\lp M^{\tfrac{4}{3}} + M^2 \rp T}. 
\end{align}
Next, to determine $c_2$ we integrate \eqref{Ohm_conv6} on the interval $[0,T]$, apply \eqref{Ohm_conv_rate1}, and obtain 

\begin{align}\label{Ohm_conv_rate3}
    &\int_0^T 
    \nu \normHs{1}{\grad \tu}^2 +
    \mu \normHs{1}{\grad \tu}^2 +
    \kappa \normHs{1}{\grad \tu}^2 \, dt 
    \\ &\leq \notag 
    C M^{2+2\beta}\lp e^{C\lp M^{\tfrac{4}{3}} + M^2 \rp T} + T \rp \epsilon^{2\gamma}.
\end{align}

Thus \eqref{thm_conv_rate2} is valid with 
\begin{align}\label{c_2}
    c_2 = C M^{2+2\beta}\lp e^{C\lp M^{\tfrac{4}{3}} + M^2 \rp T} + T \rp.
\end{align}
From \eqref{Ohm_conv_rate1} and \eqref{Ohm_conv_rate3} we conclude that $u, b $ converge to $U, B$, respectively, 
in $L^\infty(0,T; V) \cap L^2(0,T; H^2\cap V) $ and that $\theta$ converges to 
$\Theta$ in $L^\infty(0,T; L^2) \cap L^2(0,T; H^1)$. This completes the proof of convergence.  \qedsymbol

\section*{Concluding Remarks}
In this manuscript, we have proposed the first known globally well-posed approximation of the Ohmic-heating system \eqref{Sys2}, and given a rigorous proof of this.  Moreover, we prove, under the assumption that the original equation has a smooth solution (the existence of which remains a challenging open problem), that strong solutions of this approximate system \eqref{Sys1} converge to strong solutions of the original Ohmic-heating system \eqref{Sys2} as the calming parameter $\epsilon$ goes to zero on the time interval of existence and uniqueness of smooth solutions to \eqref{Sys2}.

We also note that one could extend this work to the $3$D case by using the idea of applying a calming function to the Ohmic heating term, with additional turbulence models or regularizations to handle the Navier-Stokes part. For instance, this could be done using an $\alpha$-model, the Voigt regularization, higher-order fractional diffusion, or even calming the advective nonlinearity; i.e., replacing the term $u\cdot\nabla u$ by $\bze(u)\cdot\nabla u$ (or $(\nabla\times u)\times\bze(u)$) as was done in the case of the 3D Navier-Stokes equations in \cite{Enlow_Larios_Wu_2023_NSE}.

While finalizing this manuscript, we realized that one could also consider a slightly different modification of the Ohmic heating term using a so-called ``soft-plus'' function in place of the absolute value; namely,
\[
\alpha\mu \bze(\sqrt{\epsilon+|\nabla\times b|^2})\sqrt{\epsilon+|\nabla\times b|^2}
\]
which would slightly improve the smoothness of this term, and would therefore improve the smoothness of $\partial_t\theta$  (see Remark \ref{rmk_reg_problem}).  We plan to explore this idea in future work.

\section*{Acknowledgments}
 \noindent
  The research of A.L. and M.E. was supported in part by NSF grants DMS-2206762 and CMMI-1953346.  
  A.L. was also supported in part by USGS Grant No. G23AC00156-01.

\end{document}